\documentclass[11pt, a4paper]{article}
\usepackage{amsmath} \usepackage{euscript}
\usepackage{mymatrix}
\usepackage{amssymb}
\begin{document}

\newcounter{bnomer} \newcounter{snomer}
\newcounter{bsnomer}
\setcounter{bnomer}{0}
\renewcommand{\thesnomer}{\thebnomer.\arabic{snomer}}
\renewcommand{\thebsnomer}{\thebnomer.\arabic{bsnomer}}
\renewcommand{\refname}{\begin{center}\large{\textbf{References}}\end{center}}

\newcommand{\sect}[1]{%
\setcounter{snomer}{0}\setcounter{bsnomer}{0}
\refstepcounter{bnomer}
\par\bigskip\begin{center}\large{\textbf{\arabic{bnomer}. {#1}}}\end{center}}
\newcommand{\sst}{%
\refstepcounter{bsnomer}
\par\bigskip\textbf{\arabic{bnomer}.\arabic{bsnomer}. }}
\newcommand{\defi}[1]{%
\refstepcounter{snomer}
\par\textbf{Definition \arabic{bnomer}.\arabic{snomer}. }{#1}\par}
\newcommand{\theo}[2]{%
\refstepcounter{snomer}
\par\textbf{Theorem \arabic{bnomer}.\arabic{snomer}. }{#2} {\emph{#1}}\hspace{\fill}$\square$\par}
\newcommand{\mtheo}[1]{%
\refstepcounter{snomer}
\par\textbf{Theorem \arabic{bnomer}.\arabic{snomer}. }{\emph{#1}}\par}
\newcommand{\theobp}[2]{%
\refstepcounter{snomer}
\par\textbf{Theorem \arabic{bnomer}.\arabic{snomer}. }{#2} {\emph{#1}}\par}
\newcommand{\theop}[2]{%
\refstepcounter{snomer}
\par\textbf{Theorem \arabic{bnomer}.\arabic{snomer}. }{\emph{#1}}
\par\textsc{Proof}. {#2}\hspace{\fill}$\square$\par}
\newcommand{\theosp}[2]{%
\refstepcounter{snomer}
\par\textbf{Theorem \arabic{bnomer}.\arabic{snomer}. }{\emph{#1}}
\par\textbf{Sketch of the proof}. {#2}\hspace{\fill}$\square$\par}
\newcommand{\exam}[1]{%
\refstepcounter{snomer}
\par\textbf{Example \arabic{bnomer}.\arabic{snomer}. }{#1}\par}
\newcommand{\deno}[1]{%
\refstepcounter{snomer}
\par\textbf{Definition \arabic{bnomer}.\arabic{snomer}. }{#1}\par}
\newcommand{\post}[1]{%
\refstepcounter{snomer}
\par\textbf{Proposition \arabic{bnomer}.\arabic{snomer}. }{#1}\hspace{\fill}$\square$\par}
\newcommand{\postp}[2]{%
\refstepcounter{snomer}
\par\textbf{Proposition \arabic{bnomer}.\arabic{snomer}. }{\emph{#1}}
\par\textsc{Proof}. {#2}\hspace{\fill}$\square$\par}
\newcommand{\lemm}[1]{%
\refstepcounter{snomer}
\par\textbf{Lemma \arabic{bnomer}.\arabic{snomer}. }{\emph{#1}}\hspace{\fill}$\square$\par}
\newcommand{\lemmp}[2]{%
\refstepcounter{snomer}
\par\textbf{Lemma \arabic{bnomer}.\arabic{snomer}. }{\emph{#1}}
\par\textsc{Proof}. {#2}\hspace{\fill}$\square$\par}
\newcommand{\coro}[1]{%
\refstepcounter{snomer}
\par\textbf{Corollary \arabic{bnomer}.\arabic{snomer}. }{\emph{#1}}\hspace{\fill}$\square$\par}
\newcommand{\mcoro}[1]{%
\refstepcounter{snomer}
\par\textbf{Corollary \arabic{bnomer}.\arabic{snomer}. }{\emph{#1}}\par}
\newcommand{\corop}[2]{%
\refstepcounter{snomer}
\par\textbf{Corollary \arabic{bnomer}.\arabic{snomer}. }{\emph{#1}}
\par\textsc{Proof}. {#2}\hspace{\fill}$\square$\par}
\newcommand{\nota}[1]{%
\refstepcounter{snomer}
\par\textbf{Remark \arabic{bnomer}.\arabic{snomer}. }{#1}\par}
\newcommand{\propp}[2]{%
\refstepcounter{snomer}
\par\textbf{Proposition \arabic{bnomer}.\arabic{snomer}. }{\emph{#1}}
\par\textsc{Proof}. {#2}\hspace{\fill}$\square$\par}
\newcommand{\hypo}[1]{%
\refstepcounter{snomer}
\par\textbf{Conjecture \arabic{bnomer}.\arabic{snomer}. }{\emph{#1}}\par}

\newcommand{\Ind}[3]{%
\mathrm{Ind}_{#1}^{#2}{#3}}
\newcommand{\Res}[3]{%
\mathrm{Res}_{#1}^{#2}{#3}}
\newcommand{\epsi}{\varepsilon}
\newcommand{\Supp}[1]{%
\mathrm{Supp}(#1)}

\newcommand{\reg}{\mathrm{reg}}
\newcommand{\sreg}{\mathrm{sreg}}
\newcommand{\codim}{\mathrm{codim}\,}
\newcommand{\chara}{\mathrm{char}\,}
\newcommand{\rk}{\mathrm{rk}\,}
\newcommand{\id}{\mathrm{id}}
\newcommand{\col}{\mathrm{col}}
\newcommand{\row}{\mathrm{row}}
\newcommand{\low}{\mathrm{low}}
\newcommand{\pho}{\hphantom{\quad}\vphantom{\mid}}
\newcommand{\wt}{\widetilde}
\newcommand{\wh}{\widehat}
\newcommand{\ad}[1]{\mathrm{ad}_{#1}}
\newcommand{\tr}{\mathrm{tr}\,}
\newcommand{\GL}{\mathrm{GL}}
\newcommand{\Mat}{\mathrm{Mat}}

\newcommand{\vfi}{\varphi}
\newcommand{\teta}{\vartheta}
\newcommand{\lee}{\leqslant}
\newcommand{\gee}{\geqslant}
\newcommand{\Fp}{\mathbb{F}}
\newcommand{\Rp}{\mathbb{R}}
\newcommand{\Zp}{\mathbb{Z}}
\newcommand{\Cp}{\mathbb{C}}
\newcommand{\ut}{\mathfrak{u}}
\newcommand{\at}{\mathfrak{a}}
\newcommand{\nt}{\mathfrak{n}}
\newcommand{\rt}{\mathfrak{r}}
\newcommand{\rad}{\mathfrak{rad}}
\newcommand{\bt}{\mathfrak{b}}
\newcommand{\gt}{\mathfrak{g}}
\newcommand{\vt}{\mathfrak{v}}
\newcommand{\pt}{\mathfrak{p}}
\newcommand{\Po}{\EuScript{P}}
\newcommand{\Uo}{\EuScript{U}}
\newcommand{\Fo}{\EuScript{F}}
\newcommand{\Do}{\EuScript{D}}
\newcommand{\Eo}{\EuScript{E}}
\newcommand{\Mo}{\mathcal{M}}
\newcommand{\Nu}{\mathcal{N}}
\newcommand{\Ro}{\mathcal{R}}
\newcommand{\Co}{\mathcal{C}}
\newcommand{\Lo}{\mathcal{L}}
\newcommand{\Ou}{\mathcal{O}}
\newcommand{\Au}{\mathcal{A}}
\newcommand{\Vu}{\mathcal{V}}
\newcommand{\Bu}{\mathcal{B}}
\newcommand{\Sy}{\mathcal{Z}}
\newcommand{\Sb}{\mathcal{F}}
\newcommand{\Gr}{\mathcal{G}}

\author{Mikhail V. Ignatyev\thanks{The research was partially supported by RFBR grant no. 11--01--90703-mob\_st}}
\date{\small Department of Algebra and Geometry\\
Samara State University\\
Samara, 443011, Ak. Pavlova, 1, Russia,\\
\texttt{mihail.ignatev@gmail.com}}
\title{\Large{Combinatorics of $B$-orbits and Bruhat--Chevalley order on involutions}} \maketitle

\vspace{-1cm}\sect{Introduction and main results}

\sst Let $S_n$ be the symmetric group on $n$ letters. Bruhat--Chevalley order on $S_n$ is fundamental in a multitude of contexts. For example, it describes the incidences among the closures of double cosets in the Bruhat decomposition of the general linear group $\GL_n(\Cp)$. An interesting subposet of Bruhat--Chevalley order is induced by the involutions, i.e., the elements of order 2 of $S_n$ (we denote this subposet by $S_n^2$). Activity around~$S_n^2$ was initiated by R. Richardson and T. Springer \cite{RichardsonSpringer}, who proved that the inverse Bruhat--Chevalley order on~$S_{2n+1}^2$ encodes the incidences among the closed orbits under the action of the Borel subgroup on the symmetric variety $\mathrm{SL}_{2n+1}(\Cp)/\mathrm{SO}_{2n+1}(\Cp)$.

The poset of involutions was also studied by F. Incitti \cite{Incitti1}, \cite{Incitti2} from a purely combinatorial point of view. In particular, he proved that this poset is graded, calculated the rank function and described the covering relations. In \cite{BagnoCherniavsky}, E. Bagno and Y. Chernavsky present another geometrical interpretation of the poset $S_n^2$, considering the action of standard Borel subgroup $B$ (i.e., the group of upper-triangular invertible matrices) of $\GL_n(\Cp)$ on symmetric matrices by congruence. Note that all geometric interpretations deal with the closures of orbits for various actions of the Borel subgroup. The purpose of the paper is to incorporate \emph{coadjoint} orbits into the picture.

Let $\nt$ be the space of strictly upper-triangular matrices and $\nt^*$ its dual space. Since $B$ acts on $\nt$ by conjugation, one can consider the dual action of $B$ on $\nt^*$. To~each involution $\sigma\in S_n^2$ one can assign the $B$-orbit $\Omega_{\sigma}\subset\nt^*$ (see Subsection \ref{sst:repr_theo_approach} for precise definitions). Our main result is as follows.

\medskip \mtheo{Let $\sigma,\tau\in S_n^2$. The orbit $\Omega_{\tau}$ is contained in the Zariski closure of $\Omega_{\sigma}$ if and only if $\tau\leq\sigma$ with respect to Bruhat--Chevalley order.\label{mtheo_0}}

\medskip Note that in \cite{Melnikov1}, \cite{Melnikov2}, \cite{Melnikov3} A. Melnikov described the incidences among the closures of\linebreak $B$-orbits on the variety of upper-triangular $2$-nilpotent matrices in combinatorial terms of so-called link patterns and rook placements. (In \cite{BoosReineke}, M. Boos and M. Reineke generalize the results of Melnikov to all 2-nilpotent matrices; see also B. Rothbach's paper \cite{Rothbach}.) In some sense, our results are ``dual'' to Melnikov's results.

\medskip The paper is organized as follows. In the rest of this Section, we define orbit $\Omega_{\sigma}$ associated to involution $\sigma$ from the perspective of representation theory, combinatorics and geometry. Namely, in Subsection \ref{sst:repr_theo_approach}, we give precise definitions and explain the role of orbits $\Omega_{\sigma}$ in A.A. Kirillov's orbit method in representation theory of the unipotent radical of $B$. In Subsection~\ref{sst:our_order}, we briefly recall Melnikov's results and define the partial order $\leq^*$ on~$S_n^2$ in combinatorial terms in the spirit of \cite{Melnikov2}. Then, we formulate Theorem~\ref{mtheo} claiming that $\leq^*$ encodes the incidences among the closures of $\Omega_{\sigma}$, $\sigma\in S_n^2$. In Subsection \ref{sst:Bruhat_order}, we formulate Theorem~\ref{mcoro} claiming that the restriction of Bruhat--Chevalley order to $S_n^2$ coincides with $\leq^*$. Next, in Subsection~\ref{sst:geom_conj_approach}, we present a conjectural approach based on the geometry of tangent cones to Schubert varieties.

\newpage In Section \ref{sect:proof_MT}, we prove Theorem \ref{mtheo}
(see Propositions \ref{prop:if} and \ref{prop:only_if}). In {Sub\-sec\-tion}~\ref{sst:Bruhat}, using Incitti's results, we prove Theorem \ref{mcoro}. This concludes the proof of our main result. Section \ref{sect:proofs_Lemmas} contains the proofs of technical (but important) Lemmas
used in the proof of Proposition~\ref{prop:only_if}. Finally, in
Section \ref{sect:remarks}, we present a formula for the dimension of $\Omega_{\sigma}$ (see Proposition~\ref{prop:dim_Omega}). We also formulate a conjecture about the closure of $\Omega_{\sigma}$ and check it in some particular cases (see {Sub\-sec\-tion}~\ref{sst:closure_conj}). A~short announcement of our results was made in
\cite{Ignatev}.

\medskip\textsc{Acknowledgements}. A part of this work was carried out during my stay at Moscow State University. I would like to thank Professor E.B. Vinberg for his hospitality. Financial support from RFBR (grant no. 11--01--90703-mob\_st) is gratefully acknowledged.

I~am very grateful to A.N. Panov
for useful discussions. I would like to express my gratitude to
E.Yu.~Smirnov. He was the first to explain me the role of rook placements
in combinatorics of the orbit closures. I also would like to thank A. Melnikov for suggesting the idea of Conjecture \ref{conj_closure}. I~thank two anonymous referees for useful comments on a previous version of the paper.

\sst\label{sst:repr_theo_approach} Let $G=\mathrm{GL}_n(\Cp)$ be the general linear group, $B$ its
standard Borel subgroup (etc.). Let $U\subset B$ be the \emph{unitriangular} group (i.e., the group of upper-triangular matrices with $1$'s on the diagonal). Group $B$ acts on $\nt$ by conjugation, so the dual
action of~$B$~on~$\nt^*$ is induced. For $g\in B$ and $\lambda\in\nt^*$, $g.\lambda$ is defined by
\begin{equation*}
(g.\lambda)(x)=\lambda(g^{-1}xg),~\text{for }x\in\nt.
\end{equation*}
Let $\Omega_{\lambda}$ denote the orbit of $\lambda\in\nt^*$ under this
action. Let $\Theta_{\lambda}$ denote the orbit of $\lambda$ under the action of~$U$. (Clearly, $\Theta_\lambda\subseteq\Omega_{\lambda}$.)

In 1962, Kirillov showed \cite{Kirillov1} that there is a bijection between the set $\nt^*/U$ of $U$-orbits on $\nt^*$ and the set $\widehat U$ of equivalence classes of unitary irreducible representations of $U$ in Hilbert spaces. (The proof was adapted for unipotent groups over finite fields by D. Kazhdan in \cite{Kazhdan1}.) Further, it turned out that all the principle questions about representations can be answered in terms of orbits (see \cite{Kirillov2} for the details). However, a complete description of $\nt^*/U$ is unknown and seems to be a very difficult problem.

An element $\sigma\in S_n$ satisfying $\sigma^2=\id$ is called
an \emph{ivolution}. Let $S_n^2$ be the set of
involutions of~$S_n$. To $\sigma\in S_n^2$ one can assign the
orbits of the groups $B$~and~$U$ by the following rule. Write $\sigma$ as a
product of disjoint cycles: $\sigma=(i_1,j_1)\ldots(i_t,j_t)$, where
$i_l>j_l$ for $1\leq l\leq t$ and $j_l<j_{l+1}$ for $1\leq l<t$. Denote
\begin{equation*}
\Phi=\{(i,j),1\leq j<i\leq n\}\subset\Zp\times\Zp
\end{equation*}
and put $\Supp{\sigma}=\bigcup_{l=1}^t\{(i_l,j_l)\}\subset\Phi$. Clearly, $\{e_{\alpha},\alpha\in\Phi\}$ is a basis of $\nt$. Here $e_\alpha=e_{j,i}$ for $\alpha=(i,j)\in\Phi$, where $e_{j,i}$ is the usual matrix unit. Hence one can consider the dual basis $\{e_{\alpha}^*,\alpha\in\Phi\}$ of $\nt^*$. Now, to each map $\xi\colon\Supp{\sigma}\to\Cp^{\times}\colon\alpha=(i_l,j_l)\mapsto\xi_l$ one can assign the $U$-orbit $\Theta_{\sigma,\xi}$ by putting $\Theta_{\sigma,\xi}=\Theta_{f_{\sigma,\xi}}$, where
\begin{equation*}
f_{\sigma,\xi}=\sum_{\alpha\in\Supp{\sigma}}\xi(\alpha)e_{\alpha}^*=\sum_{l=1}^t\xi_le_{j_l,i_l}^*.
\end{equation*}
(If $\sigma=\mathrm{id}$, then $\Supp{\sigma}=\emptyset$ and $f_{\sigma,\xi}=0$.) We say that $\Theta_{\sigma,\xi}$ is \emph{associated} with $\sigma$ and $\xi$. Set $\xi_0(\alpha)=1$ for all $\alpha\in\Supp{\sigma}$, and $\Omega_{\sigma}=\Omega_{f\sigma,\xi_0}$. (In other words, $f_{\sigma,\xi_0}=\sum_{\alpha\in\Supp{\sigma}}e_{\alpha}^*$.) Lemma~\ref{lemm:unip_orbits} shows that $\Omega_{\sigma}=\bigcup\Theta_{\sigma,\xi}$, where the union is taken over all maps $\xi\colon\Supp{\sigma}\to\Cp^{\times}$.

\medskip It turned out that almost all $U$-orbits on $\nt^*$ studied so far are associated with involutions.
\exam{i) Being an orbit of a connected unipotent group on an affine variety, any $U$-orbit is a Zariski-closed irreducible subvariety of $\nt^*$. Let $\Theta$ be an orbit of maximal dimension (such an orbit is called \emph{regular}). Then either $\Theta=\Theta_{w_0,\xi}$ or $\Theta=\Theta_{w_1,\xi}$ for some $\xi$ (in the last case $n$ must be even). Here $w_0=(n,1)(n-1,2)\ldots(n-n_0+1,n_0)$, $n_0=[n/2]$, and $\Supp{w_1}=\Supp{w_0}\setminus\{(n-n_0+1,n_0)\}$. Conversely, all $\Theta_{w_1,\xi}$'s are regular \cite[\S9, Example 2]{Kirillov1}.

ii) An orbit $\Theta\subset\nt^*$ is called \emph{subregular} if it has the second maximal dimension. Pick $1\leq j<n_0$ and put $\sigma$ to be the involution such that
\begin{equation*}
\Supp{\sigma}=(\Supp{w_0}\setminus\{(n-j+1,j),(n-j,j+1)\})\cup\{(n-j+1,j+1),(n-j,j)\}.
\end{equation*}
Then $\Theta_{\sigma,\xi}$ is subregular for all $\xi$. Subregular orbits were described by Panov in \cite{Panov}.

iii) Let $\alpha=(i,j)\in\Phi$. The orbit of $e_{\alpha}^*$ is called \emph{elementary}. Evidently, it is associated with the involution $\sigma=(i,j)\in S_n^2$. Elementary orbits are described in \cite{Mukherjee1}.\label{exam:coad_unip_orbits}}

Thus, orbits associated with involutions play an important role in representation theory. (See \cite{Andre1}, \cite{Andre2}, \cite{AndreNeto1}, \cite{Ignatev1} and \cite{Ignatev2} for further examples and generalizations to other unipotent algebraic groups.) They were completely described by Panov in \cite{Panov}. In particular, for a given orbit $\Theta_{\sigma,\xi}$, he presented the set of equations defining this orbit as a closed subvariety of $\nt^*$. On the contrary, $B$-orbits $\Omega_{\sigma}$ are \emph{not} closed, so the natural question arises: given two orbits $\Omega_{\tau}$ and $\Omega_{\sigma}$, $\sigma,\tau\in S_n^2$, when $\Omega_{\tau}\subseteq\overline{\Omega}_{\sigma}$? (Here $\overline{Z}$ denotes the Zariski closure of a subset $Z\subseteq\nt^*$.) By Theorem \ref{mtheo_0}, this occurs if and only if $\tau\leq_B\sigma$, where $\leq_B$ denotes Bruhat--Chevalley order.

\sst\label{sst:our_order} Let $\Nu\subset\nt$ be the variety of upper-triangular matrices of square zero:
\begin{equation*}
\Nu=\{X\in\nt\mid X^2=0\}.
\end{equation*}
Group $B$ acts on $\Nu$ by conjugation. For a given $X\in\Nu$, let $\Ou_X$ denote the orbit of $X$ under this action. To $\sigma\in S_n^2$ one can assign the
orbit $\Ou_{\sigma}$ by the following rule. Write $\sigma$ as a
product of disjoint cycles: $\sigma=(i_1,j_1)\ldots(i_t,j_t)$, where
$i_l>j_l$ for $1\leq l\leq t$ and $j_l<j_{l+1}$ for $1\leq l<t$.
Denote by~$X_{\sigma}\in\Nu$ the matrix of the form $X_{\sigma}=\sum_{\alpha\in\Supp{\sigma}}e_{\alpha}=\sum_{l=1}^te_{j_l,i_l}$, and put $\Ou_{\sigma}=\Ou_{X_{\sigma}}$. By \cite[Theorem~2.2]{Melnikov1}, one has
\begin{equation*}
\Nu=\bigsqcup_{\sigma\in S_n^2}\Ou_{\sigma}.
\end{equation*}

To each $\sigma\in S_n^2$ one can also assign the matrix
$R_{\sigma}$ by putting
\begin{equation*}
(R_{\sigma})_{i,j}=\rk\pi_{i,j}(X_{\sigma}),
\end{equation*}
where $\pi_{i,j}\colon\mathrm{Mat}_n(\Cp)\to\mathrm{Mat}_n(\Cp)$ acts on a matrix by replacing all
entries of the first $(i-1)$ rows and the last $(n-j)$ columns by
zeroes. Let us define a partial order on $S_n^2$. Given
$\sigma,\tau\in S_n^2$, we put $\sigma\leq\tau$ if $R_{\sigma}\leq
R_{\tau}$, i.e., $(R_{\sigma})_{i,j}\leq(R_{\tau})_{i,j}$ for all
$1\leq i<j\leq n$.

\medskip\exam{Let $n=5$, $\sigma=(3,1)(5,2)$,
$\tau=(2,1)(4,3)\in S_5^2$. Then
\begin{equation*}
\postdisplaypenalty=10000
\begin{split}
&X_{\sigma}=\begin{pmatrix}0&0&1&0&0\\
0&0&0&0&1\\
0&0&0&0&0\\
0&0&0&0&0\\
0&0&0&0&0
\end{pmatrix},\text{ }
X_{\tau}=\begin{pmatrix}0&1&0&0&0\\
0&0&0&0&0\\
0&0&0&1&0\\
0&0&0&0&0\\
0&0&0&0&0
\end{pmatrix},\\
&R_{\sigma}=\begin{pmatrix}0&0&1&1&2\\
0&0&0&0&1\\
0&0&0&0&0\\
0&0&0&0&0\\
0&0&0&0&0
\end{pmatrix},\text{ }
R_{\tau}=\begin{pmatrix}0&1&1&2&2\\
0&0&0&1&1\\
0&0&0&1&1\\
0&0&0&0&0\\
0&0&0&0&0
\end{pmatrix},
\end{split}
\end{equation*}
so $R_{\sigma}\leq R_{\tau}$ and $\sigma\leq\tau$.} \nota{Note that
this partial order has an interpretation in terms of so-called
\emph{rook {place\-ments}}. Namely, $X_{\sigma}$ can be treated as a
rook placement on the triangle board with boxes labeled by
pairs $(i, j)$, $1\leq i<j\leq n$: by definition, there is a rook in
the $(i, j)$th box if and only if $(X_{\sigma})_{i,j}=1$. Then
$(R_{\sigma})_{i,j}$ is just the number of rooks located
non-strictly to the South-West of the $(i,j)$th box.}

\medskip As above, let $\overline{Z}$ be the closure of a subset
$Z\subseteq\mathrm{Mat}_n(\Cp)$ with respect to Zariski topology. By \cite[Theorem~3.5]{Melnikov2}, one has the following nice combinatorial
description of the orbit closures in $\Nu$:
\begin{equation}
\Ou_{\tau}\subseteq\overline{\Ou}_{\sigma}\text{ if and only if
}\tau\leq\sigma.\label{formula:Melnikov_result}
\end{equation}
In \cite{Melnikov3}, an interpretation of this result in terms of
link patterns is given.

\medskip Now, let $\nt_-$ be the space of strictly lower-triangular
matrices (with zeroes on the diagonal). We can identify it with $\nt^*$ by putting
\begin{equation*}
\lambda(x)=\langle\lambda,x\rangle=\tr{\lambda x},\text{
}\lambda\in\nt_-,\text{ }x\in\nt.
\end{equation*}
Thus, in the sequel we identify $\nt^*$ with $\nt_-$. Note that under this identification, $e_{\alpha}^*=e_{i,j}$ for all $\alpha=(i,j)\in\Phi$, and $\Omega_{\lambda}=\{(g\lambda g^{-1})_{\mathrm{low}},g\in B\}$, where $A_{\mathrm{low}}$ denotes the strictly lower-triangular part of $A$, that is
\begin{equation*}
(A_{\low})_{i,j}=\begin{cases}A_{i,j},&\text{if }i>j,\\
0&\text{otherwise.}
\end{cases}
\end{equation*}

Let $\sigma\in S_n^2$. Then $f_{\sigma,\xi_0}$ is identified with $X_{\sigma}^t$, so $\Omega_{\sigma}$ is identified with $\Omega_{X_{\sigma}^t}$,
where $X_{\sigma}^t\in\nt^*$ denotes the transposed matrix to
$X_{\sigma}$.
In fact, our goal is to describe $\overline{\Omega}_{\sigma}$ in combinatorial terms. To do
this, let us define another partial order on $S_n^2$. Given
$\sigma$, $\tau\in S_n^2$, put $\sigma\leq^*\tau$ if
$R_{\sigma}^*\leq R_{\tau}^*$, i.e., $(R_{\sigma}^*)_{i,j}\leq
(R_{\tau}^*)_{i,j}$ for all $1\leq j<i\leq n$. Here
$R_{\sigma}^*\in\nt^*$ is the matrix defined by the rule
\begin{equation*}
(R_{\sigma}^*)_{i,j}=\begin{cases}\rk\pi_{i,j}(X_{\sigma}^t),&\text{if
}1\leq j<i\leq n,\\
0&\text{otherwise}.
\end{cases}
\end{equation*}
As above, $\pi_{i,j}\colon\mathrm{Mat}_n(\Cp)\to\mathrm{Mat}_n(\Cp)$ acts on a matrix by replacing
all entries of the first $(i-1)$ rows and the last $(n-j)$ columns
by zeroes.

\exam{Let $n=5$, $\sigma=(4,1)(5,2)$, $\tau=(5,1)(4,2)\in S_5^2$.
Then
\begin{equation*}
\begin{split}
&X_{\sigma}^t=\begin{pmatrix}0&0&0&0&0\\
0&0&0&0&0\\
0&0&0&0&0\\
1&0&0&0&0\\
0&1&0&0&0
\end{pmatrix},\text{ }
X_{\tau}^t=\begin{pmatrix}0&0&0&0&0\\
0&0&0&0&0\\
0&0&0&0&0\\
0&1&0&0&0\\
1&0&0&0&0
\end{pmatrix},\\
&R_{\sigma}^*=\begin{pmatrix}0&0&0&0&0\\
1&0&0&0&0\\
1&2&0&0&0\\
1&2&2&0&0\\
0&1&1&1&0
\end{pmatrix},\text{ }
R_{\tau}^*=\begin{pmatrix}0&0&0&0&0\\
1&0&0&0&0\\
1&2&0&0&0\\
1&2&2&0&0\\
1&1&1&1&0
\end{pmatrix},
\end{split}
\end{equation*}
so $R_{\sigma}^*\leq R_{\tau}^*$ and $\sigma\leq^*\tau$.}\nota{Of
course, this partial order has an interpretation in terms of
rook placements. Namely, $X_{\sigma}^t$ can be treated as a rook
placement on the triangle board with boxes labeled by pairs
$(i, j)$, $1\leq j<i\leq n$: by definition, there is a rook in the
$(i, j)$th box if and only if $(X_{\sigma}^t)_{i,j}=1$. Then
$(R_{\sigma}^*)_{i,j}$, $i>j$, is just the number of rooks located
non-strictly to the South-West of the $(i,j)$th box.}

\medskip The following theorem plays a key role in the proof of the main result of the paper (cf. (\ref{formula:Melnikov_result})).
\mtheo{Let $\sigma,\tau$ be involutions in $S_n$ and
$\Omega_{\sigma},\Omega_{\tau}$ the corresponding $B$-orbits in
$\nt^*$. Then\label{mtheo}
\begin{equation*}
\Omega_{\tau}\subseteq\overline{\Omega}_{\sigma}\text{ if and only
if }\tau\leq^*\sigma.
\end{equation*}}
The proof will be presented in the next Section (see Proposition \ref{prop:if} for the proof of ``only if'' direction and Proposition \ref{prop:only_if} for the proof of ``if'' direction).

\nota{Note that there is \emph{no} natural analogue of the variety $\Nu$ in the space $\nt^*$. Actually, one can put
\begin{equation*}
\Nu^*=\bigsqcup_{\sigma\in S_n^2}\Omega_{\sigma}.
\end{equation*}
Clearly, this subset of $\nt^*$ is stable under the action of $B$, but it is neither open nor closed, if $n>2$. (For $n=2$, $\Nu=\nt^*$.) Indeed, it contains the orbit $\Omega_{w_0}$, where $w_0=(n,1)(n-1,2)\ldots(n-n_0+1,n_0)$, $n_0=[n/2]$ (as in Example~\ref{exam:coad_unip_orbits}i)). It~follows from \cite[\S9, Example 2]{Kirillov1} and Lemma \ref{lemm:unip_orbits} that $y\in\nt^*$ belongs to $\Omega_{w_0}$ if and only if $\Delta_i(y)\neq0$ for all $1\leq i\leq n_0$. Here
\begin{equation*}
\Delta_i(y)=\begin{vmatrix}y_{n-i+1,1}&y_{n-i+1,2}&\ldots&y_{n-i+1,i}\\
y_{n-i+2,1}&y_{n-i+2}&\ldots&y_{n-i+2,i}\\
\vdots&\vdots&\ddots&\vdots\\
y_{n,1}&y_{n,2}&\ldots&y_{n,i}
\end{vmatrix}.
\end{equation*} Hence $\Omega_{w_0}$ is an open subset of $\nt^*$, so $\overline{\Nu^*}=\overline{\Omega}_{w_0}=\nt^*$ and $\Nu^*$ is not closed.

On the other hand, suppose $\Nu^*$ is open. Consider $V=\{y\in\nt^*\mid y_{i,j}=0,\text{ if } i>3\}$. Then $\Nu^*\cap V$ must be an open subset of $V$. However, Lemma \ref{lemm:unip_orbits}  together with \cite[Theorem 1.4]{Panov} imply that $$\Nu^*\cap V=\{y\in V\mid y_{3,1}\neq0\}\cup\{y\in V\mid y_{2,1}=y_{3,1}=0\}\cup\{y\in V\mid y_{3,1}=y_{3,2}=0\},$$ which is obviously not an open subset of $V$, a contradiction. Note, however, that $\Nu^*$ is an irreducible constructive subset of $\nt^*$ (as a union of orbits containing a dense subset of $\nt^*$). Note also that, unlike of the adjoint case considered by Melnikov, the closure of a given $\Omega_{\sigma}$, $\sigma\in S_n^2$, is \emph{not} a subset of~$\Nu^*$ (see Subsection \ref{sst:closure_conj} for a conjectural description of $\overline{\Omega}_{\sigma}$).}

\sst\label{sst:Bruhat_order} Recall that the \emph{Bruhat}--\emph{Chevalley order} $\leq_B$ on $S_n$ is defined in terms of the inclusion {re\-la\-tion\-ships} of double cosets in $\GL_n(\Cp)$. Namely, $G=\GL_n(\Cp)=\bigcup_{w\in S_n}B\dot wB$, where $\dot w$ denotes the permutation matrix corresponding to $w$. Let $v,w\in S_n$. By definition, $v\leq_Bw$ if $B\dot vB\subseteq\overline{B\dot wB}$. Let $w=s_1\ldots s_l$ be a \emph{reduced} expression of $w$ as a product of \emph{simple reflections} $s_i=(i,i+1)\in S_n$, $1\leq i\leq n-1$, and $l(w)=l$. It's well-known that $$\{v\in S_n\mid v\leq_B w\}=\{s_{i_1}\ldots s_{i_k},\text{ }1\leq i_1<\ldots<i_k\leq l\}.$$

Further, let $X\in\mathrm{Mat}_n(\Cp)$ be an arbitrary $0$--$1$ matrix with at most one $1$ in every row and every column. Denote by $R(X)$ the matrix such that
\begin{equation*}
R(X)_{i,j}=\rk\pi_{i,j}(X),\text{ }1\leq i,j\leq n
\end{equation*}
(see the previous Subsection for the definition of $\pi_{i,j}$).
\nota{Notice that $R(X)_{i,j}$ is just the number of rooks located
non-strictly to the South-West of the $(i,j)$th box. As above, for a given matrix $A\in\mathrm{Mat}_n(\Cp)$, let $A_{\mathrm{low}}$ denote the strictly {lo\-wer-tri\-an\-gular} part of~$A$. Then $R_{\sigma}^*=R(\dot\sigma)_{\mathrm{low}}$.\label{nota:rook_placements}}

Let $v$, $w\in S_n$. Then (see, e.g., \cite{Proctor})
\begin{equation*}
v\leq_B w\text{ if and only if }R(\dot v)\leq R(\dot w),\text{i.e., } R(\dot v)_{i,j}\leq R(\dot w)_{i,j}
\text{ for all }1\leq i,j\leq n.\label{formula:Bruhat}
\end{equation*}
Suppose $\sigma$, $\tau\in S_n^2$. It follows immediately from Remark~\ref{nota:rook_placements} that $\tau\leq_B\sigma$ implies $\tau\leq^*\sigma$. In fact, the second ingredient of the proof of Theorem \ref{mtheo_0} is the fact that these conditions are equivalent, i.e., the order on $S_n^2$ induced by Bruhat--Chevalley order coincides with~$\leq^*$.\mtheo{Let $\sigma,\tau$ be involutions in $S_n$. Then\label{mcoro}
\begin{equation*}
\tau\leq^*\sigma\text{ if and only if }\tau\leq_B\sigma.
\end{equation*}}
The proof based on the computing the covering relations for $\leq^*$ and on Incitti's results is presented  in Subsection \ref{sst:Bruhat}. Note that this Theorem together with Theorem \ref{mtheo} imply our main result.

\sst\label{sst:geom_conj_approach} Before starting with the proof of Theorem \ref{mtheo_0}, we will briefly describe another (conjectural) approach to orbits
associated with involutions in terms of tangent cones to Schubert varieties. Since\linebreak $G=\bigcup_{w\in S_n}B\dot w B$, the \emph{flag variety} $\Fo=G/B$
can be decomposed into the union $\Fo=\bigcup_{w\in S_n}X_w^{\circ}$, where $X_w^{\circ}=B\dot wB/B$
is called the \emph{Schubert cell}. By definition, the \emph{Schubert variety} $X_w$ is the closure of
$X_w^{\circ}$ in $\Fo$ with respect to Zariski topology. Note that $p=X_{\id}=B/B$ is contained in $X_w$ for all $w\in S_n$. One has $X_w\subseteq X_{w'}$ if and only if $w\leq_Bw'$. Let $T_w$ be the tangent space and $C_w$ the tangent cone to $X_w$ at
the point $p$ (see \cite{BileyLakshmibai} for detailed constructions); by definition, $C_w\subseteq T_w$
and if $p$ is a~regular point of $X_w$, then $C_w=T_w$. Of course, if $w\leq_Bw'$, then $C_w\subseteq C_{w'}$.

Let $T=T_p\Fo$ be the tangent space to $\Fo$ at $p$. It can be naturally identified with $\nt^*$ in the following way. Since $\Fo=G/B$, $T$ is isomorphic to the factor $\gt/\bt$, where $\gt=\mathrm{Mat}_n(\Cp)$ is the Lie algebra of $G$ and $\bt=\langle e_{i,j},1\leq i\leq j\leq n\rangle_{\Cp}$ is the Lie algebra of $B$. In turn, $\gt/\bt$ is naturally isomorphic to $\nt_-=\nt^*$. Next, $B$ acts on $\Fo$ by left multiplications. Since $p$ is invariant under this action, the action on $T=\nt^*$ is induced. One can easily check that this action coincides with the action of $B$ on $\nt^*$ defined above \cite[Section 3, Theorem 1]{Kirillov3}. Further, the tangent cone $C_w\subseteq T_w\subseteq T=\nt^*$ is $B$-invariant, so it splits into a union of $B$-orbits.

It is well-known that $C_w$ is a subvariety of $T_w$ of dimension $\dim C_w=l(w)$ \cite[Chapter 2, Section~2.6]{BileyLakshmibai}. Let $\sigma\in S_n^2$. $\overline{\Omega}_{\sigma}$ is irreducible as the closure of an orbit. By Proposition \ref{prop:dim_Omega}, $\dim\overline{\Omega}_{\sigma}=\dim\Omega_{\sigma}=l(\sigma)$, so $\overline{\Omega}_{\sigma}$ is an irreducible component of $C_{\sigma}$ of maximal dimension. For $n\leq5$, $C_{\sigma}=\overline{\Omega}_{\sigma}$ for all $\sigma\in S_n^2$. (See \cite{PanovEliseev} for an explicit description of tangent cones.) Unfortunately, we can not prove the irreducibility of $C_{\sigma}$ for all $\sigma\in S_n^2$ for an arbitrary~$n$. On the other hand, we do not know counterexamples to the equality $C_{\sigma}=\overline{\Omega}_{\sigma}$. This allows us to formulate

\medskip\hypo{Let $\sigma\in S_n$ be an involution. Then the closure of the $B$-orbit $\Omega_{\sigma}\subset\nt^*$ coincides with the tangent cone $C_w$ to the Schubert variety $X_w$ at the point $p=B/B$.\label{conj_cones}}

\medskip Note that this conjecture implies that if $\tau\leq_B\sigma$, then $\Omega_{\tau}\subseteq\overline{\Omega}_{\sigma}$.


\sect{Proof of the Main Theorem}\label{sect:proof_MT}
\sst\label{sst:Near} The goal of this Subsection is to prove the
``only if'' direction of Theorem \ref{mtheo}. Fix an involution
$\sigma\in S_n^2$. Recall notation from Subsection~\ref{sst:repr_theo_approach}. Let $D\subset B$ be the subgroup of diagonal matrices. Clearly, $B=U\rtimes D$. (In other words, for a given $g\in B$,
there exist unique $u\in U$, $d\in D$ such that $g=ud$.)

\medskip\lemmp{One
has\footnote{Cf. \cite[Subsection
3.3]{Melnikov1}.} $\Omega_{\sigma}=\bigcup_{\xi\colon
\Supp{\sigma}\to\Cp^{\times}}\Theta_{\sigma,\xi}$.\label{lemm:unip_orbits}}{Let
$\xi\colon\Supp{\sigma}\to\Cp^{\times}$ be a map. If
$d=1_n+\sum_{l=1}^t(\xi_l-1)e_{i_l,i_l}\in D$, then
$d.X_{\sigma}^t=f_{\sigma,\xi}$, so $\Theta_{\sigma,\xi}\subset\Omega_{\sigma}$.
On the other hand, let $g$ be an element of $B$. Then there exist
$u\in U$, $d\in D$ such that $g=ud$, so $g.X_{\sigma}^t=u.f_{\sigma,\xi}$,
where $\xi(i_l,j_l)=g_{i_l,i_l}/g_{j_l,j_l}$. Thus,
$g.X_{\sigma}^t\in\Theta_{\sigma,\xi}$.}\medskip\lemmp{Let
$\lambda\in\Omega_{\sigma}$. Then $\rk\pi_{i,j}(\lambda)
=(R_{\sigma}^*)_{i,j}$ for all $1\leq j<i\leq
n$.\label{lemm:rank_through_orbit}}{Fix a map
$\xi\colon\Supp{\sigma}\to\Cp^{\times}$. Lemma \ref{lemm:unip_orbits} shows
that it's enough to check that if $u\in U$, $f\in\nt^*$, then
$\rk\pi_{i,j}(u.f)=\rk\pi_{i,j}(f)$ for all $1\leq j<i\leq
n$, because
$\rk\pi_{i,j}(f_{\sigma,\xi})=\rk\pi_{i,j}(X_{\sigma}^t)=(R_{\sigma}^*)_{i,j}$.
Pick an element $u\in U$. It's well-known that there exist
$\alpha_{j,i}\in\Cp$ such that
\begin{equation*}
u=\prod_{(i,j)\in\Phi}x_{j,i}(\alpha_{j,i}),
\end{equation*}
where $x_{j,i}(\alpha_{j,i})=1_n+\alpha_{j,i}e_{j,i}$ (the product
is taken in any fixed order). Hence we can assume
$u=x_{j,i}(\alpha)$ for some $(i, j)\in\Phi$, $\alpha\in\Cp$. Then
\begin{equation*}
(u.f)_{r,s}=\begin{cases}f_{j,s}+\alpha f_{i, s},&\text{if
}r=j\text{ and }1\leq s<j,\\
f_{r,i}-\alpha f_{r, j},&\text{if
}s=i\text{ and }i<r\leq n,\\
f_{r,s}&\text{otherwise}.
\end{cases}
\end{equation*}
Hence if $r>j$ and $s<i$, then $\pi_{r,s}(u.f)=\pi_{r,s}(f)$. If
$r\leq j$ (and so $s<r\leq j<i$), then the $j$th row of
$\pi_{r,s}(u.f)$ is obtained from the $j$th row of
$\pi_{r,s}(f)$ by adding the $i$th row of $\pi_{r,s}(f)$
multiplied by $\alpha$. Similarly, if $s\geq i$ (and so $r>s\geq
i>j$), then the $i$th column of $\pi_{r,s}(u.f)$ is obtained from
the $i$th column of $\pi_{r,s}(f)$ by subtracting the $j$th
column of $\pi_{r,s}(f)$ multiplied by $\alpha$. In both cases,
$\rk\pi_{r,s}(u.f)=\rk\pi_{r,s}(f)$, as required.}

\medskip\propp{Let $\sigma$, $\tau$ be involutions in $S_n$. If
$\Omega_{\tau}\subseteq\overline{\Omega}_{\sigma}$, then\footnote{Cf. \cite[Lemma
3.6]{Melnikov2}.}
$\tau\leq^*\sigma$.\label{prop:if}}{Suppose
$\sigma\ngeq^*\tau$, This means that there exists
$(i,j)\in\Phi$ such that $(R_{\sigma}^*)_{i,j}<(R_{\tau}^*)_{i,j}$.
Denote
\begin{equation*}
Z=\{f\in\nt^*\mid\rk\pi_{r,s}(f)\leq(R_{\sigma}^*)_{r,s}\text{ for
all }(r,s)\in\Phi\}.
\end{equation*}
Clearly, $Z$ is closed with respect to Zariski topology. Lemma
\ref{lemm:rank_through_orbit} shows that $\Omega_{\sigma}\subseteq
Z$, so $\overline{\Omega}_{\sigma}\subseteq Z$. But
$X_{\tau}^t\notin Z$, hence $\Omega_{\tau}\nsubseteq Z$, a
contradiction.}

\sst Now, let us start with the proof of much more difficult
``if'' direction of Theorem~\ref{mtheo}.\linebreak First, we need some
more notation (cf. \cite[Subsections 3.7--3.14]{Melnikov2}). There
exists a natural partial order  on $\Phi$. Namely, given $(a,b)$,
$(c,d)\in\Phi$, we set $(a,b)\leq(c,d)$ if $a\leq c$ and $b\geq d$;
we also set $(a,b)>(c,d)$, if $(a,b)\geq(c,d)$ and $(a,b)\neq(c,d)$.
Let $\sigma\in S_n^2$ and $(i, j)\in\Supp{\sigma}$, i.e., $i>j$ and
$\sigma(i)=j$. Denote
\begin{equation*}
m=\min\{s\mid j<s<i\text{ and }\sigma(s)=s\}.
\end{equation*}
Suppose $m$ exists. Further, suppose that there are \emph{no}
$(p,q)\in\Supp{\sigma}$ such that $(i,j)>(p,q)$,
$(i,m)\ngtr(p,q)$. Then denote by $\sigma_{(i,j)}^{\to}\in S_n^2$
the involution such that
\begin{equation*}
\Supp{\sigma_{(i,j)}^{\to}}=(\Supp{\sigma}\setminus\{(i,j)\})\cup\{(i,m)\}.
\end{equation*}

\exam{It's very convenient to draw the corresponding $X$'s as rook
placements. For example, if $n=8$, $\sigma=(3,1)(8,2)(7,6)$, then
$\sigma_{(8,2)}^{\to}=(3,1)(8,4)(7,6)$, so
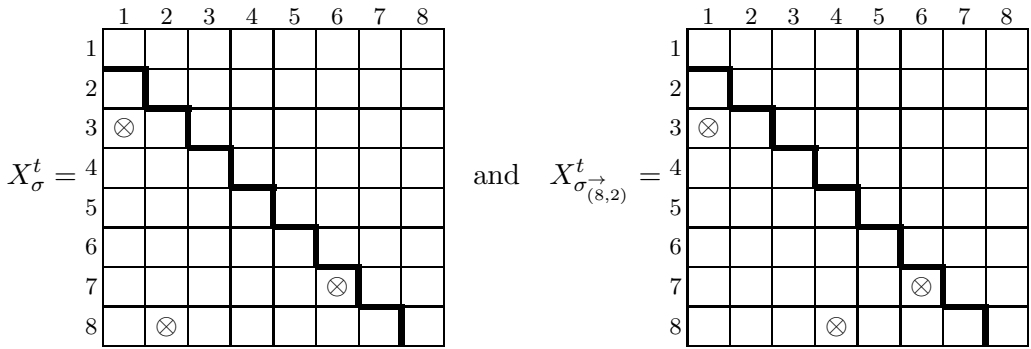
\begin{equation*}X_{\sigma}^t=
\mymatrix{ \pho& \pho& \pho& \pho& \pho& \pho& \pho& \pho\\
\Top{2pt}\Rt{2pt}\pho& \pho& \pho& \pho& \pho& \pho& \pho& \pho\\
\otimes& \Top{2pt}\Rt{2pt}\pho& \pho& \pho& \pho& \pho& \pho& \pho\\
\pho& \pho& \Top{2pt}\Rt{2pt}\pho& \pho& \pho& \pho& \pho& \pho\\
\pho& \pho& \pho& \Top{2pt}\Rt{2pt}\pho& \pho& \pho& \pho& \pho\\
\pho& \pho& \pho& \pho& \Top{2pt}\Rt{2pt}\pho& \pho& \pho& \pho\\
\pho& \pho& \pho& \pho& \pho& \Top{2pt}\Rt{2pt}\otimes& \pho& \pho\\
\pho& \otimes& \pho& \pho& \pho& \pho& \Top{2pt}\Rt{2pt}\pho& \pho\\
}\quad\text{and}\quad X_{\sigma_{(8,2)}^{\to}}^t=
\mymatrix{ \pho& \pho& \pho& \pho& \pho& \pho& \pho& \pho\\
\Top{2pt}\Rt{2pt}\pho& \pho& \pho& \pho& \pho& \pho& \pho& \pho\\
\otimes& \Top{2pt}\Rt{2pt}\pho& \pho& \pho& \pho& \pho& \pho& \pho\\
\pho& \pho& \Top{2pt}\Rt{2pt}\pho& \pho& \pho& \pho& \pho& \pho\\
\pho& \pho& \pho& \Top{2pt}\Rt{2pt}\pho& \pho& \pho& \pho& \pho\\
\pho& \pho& \pho& \pho& \Top{2pt}\Rt{2pt}\pho& \pho& \pho& \pho\\
\pho& \pho& \pho& \pho& \pho& \Top{2pt}\Rt{2pt}\otimes& \pho& \pho\\
\pho& \pho& \pho& \otimes& \pho& \pho& \Top{2pt}\Rt{2pt}\pho& \pho\\
}\
\end{equation*}

\medskip\noindent Here we denote rooks by $\otimes$'s.}

\medskip Similarly, suppose
\begin{equation*}
m=\max\{r\mid j<r<i\text{ and }\sigma(r)=r\}
\end{equation*}
exists. Further, suppose that there are no $(p,q)\in\Supp{\sigma}$
such $(i,j)>(p,q)$, $(m,j)\ngtr(p,q)$. Then
denote by $\sigma_{(i,j)}^{\uparrow}\in S_n^2$ the involution such
that
\begin{equation*}
\Supp{\sigma_{(i,j)}^{\uparrow}}=(\Supp{\sigma}\setminus\{(i,j)\})\cup\{(m,j)\}.
\end{equation*}

\exam{Let $n=8$, $\sigma=(4,1)(7,2)(8,6)$, then
$\sigma_{(7,2)}^{\uparrow}=(4,1)(5,2)(8,6)$, so
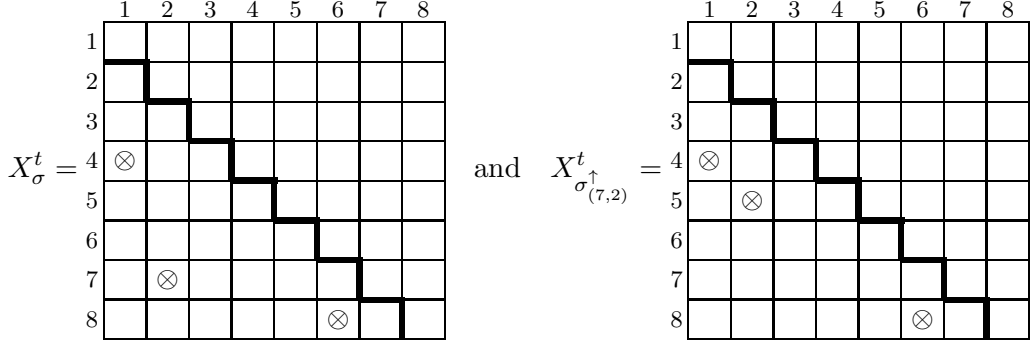
\begin{equation*}X_{\sigma}^t=
\mymatrix{ \pho& \pho& \pho& \pho& \pho& \pho& \pho& \pho\\
\Top{2pt}\Rt{2pt}\pho& \pho& \pho& \pho& \pho& \pho& \pho& \pho\\
\pho& \Top{2pt}\Rt{2pt}\pho& \pho& \pho& \pho& \pho& \pho& \pho\\
\otimes& \pho& \Top{2pt}\Rt{2pt}\pho& \pho& \pho& \pho& \pho& \pho\\
\pho& \pho& \pho& \Top{2pt}\Rt{2pt}\pho& \pho& \pho& \pho& \pho\\
\pho& \pho& \pho& \pho& \Top{2pt}\Rt{2pt}\pho& \pho& \pho& \pho\\
\pho& \otimes& \pho& \pho& \pho& \Top{2pt}\Rt{2pt}\pho& \pho& \pho\\
\pho& \pho& \pho& \pho& \pho& \otimes& \Top{2pt}\Rt{2pt}\pho& \pho\\
}\quad\text{and}\quad X_{\sigma_{(7,2)}^{\uparrow}}^t=
\mymatrix{ \pho& \pho& \pho& \pho& \pho& \pho& \pho& \pho\\
\Top{2pt}\Rt{2pt}\pho& \pho& \pho& \pho& \pho& \pho& \pho& \pho\\
\pho& \Top{2pt}\Rt{2pt}\pho& \pho& \pho& \pho& \pho& \pho& \pho\\
\otimes& \pho& \Top{2pt}\Rt{2pt}\pho& \pho& \pho& \pho& \pho& \pho\\
\pho& \otimes& \pho& \Top{2pt}\Rt{2pt}\pho& \pho& \pho& \pho& \pho\\
\pho& \pho& \pho& \pho& \Top{2pt}\Rt{2pt}\pho& \pho& \pho& \pho\\
\pho& \pho& \pho& \pho& \pho& \Top{2pt}\Rt{2pt}\pho& \pho& \pho\\
\pho& \pho& \pho& \pho& \pho& \otimes& \Top{2pt}\Rt{2pt}\pho& \pho\\
}\end{equation*}}

Now, denote by $M(\sigma)$ the set of minimal elements of
$\Supp{\sigma}$ with respect to the partial order defined above. If
$(i,j)\in M(\sigma)$, then we denote also by $\sigma_{(i,j)}^-\in
S_n^2$ the involution such that
\begin{equation*}
\Supp{\sigma_{(i,j)}^-}=\Supp{\sigma}\setminus\{(i,j)\}.
\end{equation*}

\medskip Then, denote by $A_{i,j}(\sigma)$, $(i,j)\in\Supp{\sigma}$,
the set of $(\alpha,\beta)\in\Supp{\sigma}$ such that
$j<\beta<i<\alpha$, $\sigma(r)\neq r$ for all $\beta<r<i$, and there are no $(p,q)\in\Supp{\sigma}$ such that either $j<q<\beta<p<i$ or
$\beta<q<i<p<\alpha$ (i.e., either $(i,j)>(p,q)$ and
$(\beta,j)\ngtr(p,q)$, or $(\alpha,\beta)>(p,q)$ and
$(\alpha,i)\ngtr(p,q)$). If $(\alpha,\beta)\in A_{i,j}(\sigma)$,
then denote by $a\sigma_{(i,j)}^{(\alpha,\beta)}\in S_n^2$ the
involution such that
\begin{equation*}
\Supp{a\sigma_{(i,j)}^{(\alpha,\beta)}}=
(\Supp{\sigma}\setminus\{(i,j),(\alpha,\beta)\})\cup\{(\beta,j),(\alpha,i)\}.
\end{equation*}

\exam{If $n=8$, $\sigma=(5,1)(6,2)(8,4)$, then $(8,4)\in
A_{6,2}(\sigma)$, $a\sigma_{(6,2)}^{(8,4)}=(5,1)(4,2)(8,6)$, so
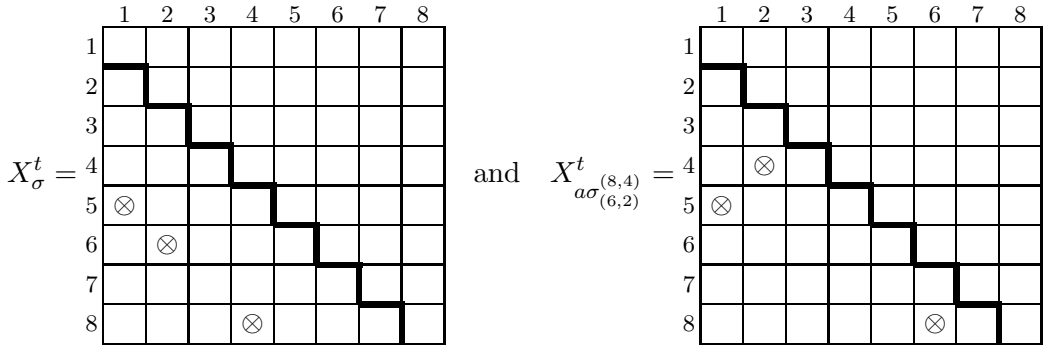
\begin{equation*}X_{\sigma}^t=
\mymatrix{ \pho& \pho& \pho& \pho& \pho& \pho& \pho& \pho\\
\Top{2pt}\Rt{2pt}\pho& \pho& \pho& \pho& \pho& \pho& \pho& \pho\\
\pho& \Top{2pt}\Rt{2pt}\pho& \pho& \pho& \pho& \pho& \pho& \pho\\
\pho& \pho& \Top{2pt}\Rt{2pt}\pho& \pho& \pho& \pho& \pho& \pho\\
\otimes& \pho& \pho& \Top{2pt}\Rt{2pt}\pho& \pho& \pho& \pho& \pho\\
\pho& \otimes& \pho& \pho& \Top{2pt}\Rt{2pt}\pho& \pho& \pho& \pho\\
\pho& \pho& \pho& \pho& \pho& \Top{2pt}\Rt{2pt}\pho& \pho& \pho\\
\pho& \pho& \pho& \otimes& \pho& \pho& \Top{2pt}\Rt{2pt}\pho& \pho\\
}\quad\text{and}\quad X_{a\sigma_{(6,2)}^{(8,4)}}^t=
\mymatrix{ \pho& \pho& \pho& \pho& \pho& \pho& \pho& \pho\\
\Top{2pt}\Rt{2pt}\pho& \pho& \pho& \pho& \pho& \pho& \pho& \pho\\
\pho& \Top{2pt}\Rt{2pt}\pho& \pho& \pho& \pho& \pho& \pho& \pho\\
\pho& \otimes& \Top{2pt}\Rt{2pt}\pho& \pho& \pho& \pho& \pho& \pho\\
\otimes& \pho& \pho& \Top{2pt}\Rt{2pt}\pho& \pho& \pho& \pho& \pho\\
\pho& \pho& \pho& \pho& \Top{2pt}\Rt{2pt}\pho& \pho& \pho& \pho\\
\pho& \pho& \pho& \pho& \pho& \Top{2pt}\Rt{2pt}\pho& \pho& \pho\\
\pho& \pho& \pho& \pho& \pho& \otimes& \Top{2pt}\Rt{2pt}\pho& \pho\\
}\end{equation*}}

\bigskip Now, denote by $B_{i,j}(\sigma)$, $(i,j)\in\Supp{\sigma}$,
the set of $(\alpha,\beta)\in\Supp{\sigma}$ such that
$(\alpha,\beta)>(i,j)$ and there are no $(p,q)\in\Supp{\sigma}$
such that $(i,j)<(p,q)<(\alpha,\beta)$. If~$(\alpha,\beta)\in
B_{i,j}(\sigma)$, then denote by
$b\sigma_{(i,j)}^{(\alpha,\beta)}\in S_n^2$ the involution such that
\begin{equation*}
\Supp{b\sigma_{(i,j)}^{(\alpha,\beta)}}=
(\Supp{\sigma}\setminus\{(i,j),(\alpha,\beta)\})\cup\{(i,\beta),(\alpha,j)\}.
\end{equation*}

\exam{Let $n=8$ and $\sigma=(8,1)(3,2)(5,4)(7,6)$. In this case, $(8,1)\in
B_{5,4}(\sigma)$ and\linebreak $b\sigma_{(5,4)}^{(8,1)}=(5,1)(3,2)(8,4)(7,6)$,
so
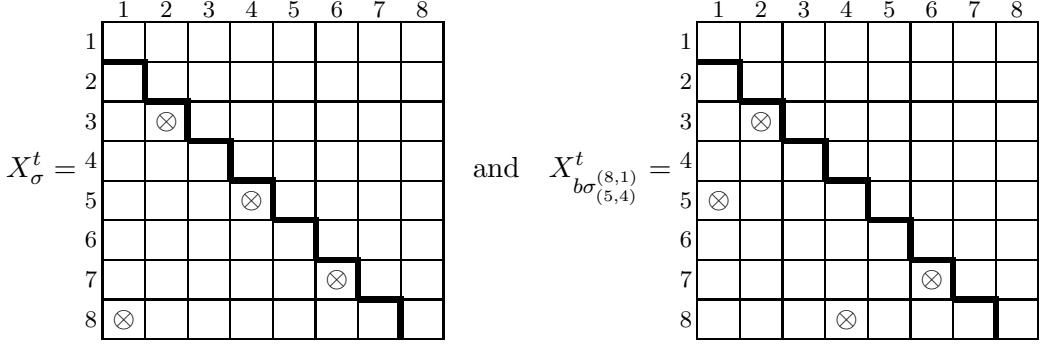
\begin{equation*}\predisplaypenalty=0
X_{\sigma}^t=
\mymatrix{ \pho& \pho& \pho& \pho& \pho& \pho& \pho& \pho\\
\Top{2pt}\Rt{2pt}\pho& \pho& \pho& \pho& \pho& \pho& \pho& \pho\\
\pho& \Top{2pt}\Rt{2pt}\otimes& \pho& \pho& \pho& \pho& \pho& \pho\\
\pho& \pho& \Top{2pt}\Rt{2pt}\pho& \pho& \pho& \pho& \pho& \pho\\
\pho& \pho& \pho& \Top{2pt}\Rt{2pt}\otimes& \pho& \pho& \pho& \pho\\
\pho& \pho& \pho& \pho& \Top{2pt}\Rt{2pt}\pho& \pho& \pho& \pho\\
\pho& \pho& \pho& \pho& \pho& \Top{2pt}\Rt{2pt}\otimes& \pho& \pho\\
\otimes& \pho& \pho& \pho& \pho& \pho& \Top{2pt}\Rt{2pt}\pho& \pho\\
}\quad\text{and}\quad X_{b\sigma_{(5,4)}^{(8,1)}}^t=
\mymatrix{ \pho& \pho& \pho& \pho& \pho& \pho& \pho& \pho\\
\Top{2pt}\Rt{2pt}\pho& \pho& \pho& \pho& \pho& \pho& \pho& \pho\\
\pho& \Top{2pt}\Rt{2pt}\otimes& \pho& \pho& \pho& \pho& \pho& \pho\\
\pho& \pho& \Top{2pt}\Rt{2pt}\pho& \pho& \pho& \pho& \pho& \pho\\
\otimes& \pho& \pho& \Top{2pt}\Rt{2pt}\pho& \pho& \pho& \pho& \pho\\
\pho& \pho& \pho& \pho& \Top{2pt}\Rt{2pt}\pho& \pho& \pho& \pho\\
\pho& \pho& \pho& \pho& \pho& \Top{2pt}\Rt{2pt}\otimes& \pho& \pho\\
\pho& \pho& \pho& \otimes& \pho& \pho& \Top{2pt}\Rt{2pt}\pho& \pho\\
}\end{equation*}}

\bigskip Finally, denote by $C_{i,j}(\sigma)$, $(i,j)\in\Supp{\sigma}$,
the set of $(\alpha,\beta)\in\Zp\times\Zp$ such that
$i>\beta>\alpha>j$, $\sigma(s)\neq s$ for all $\beta>s>\alpha$, and
if $(p,q)\in\Supp{\sigma}$, $(i,j)>(p,q)$, $(\alpha,j)\ngtr(p,q)$,
then $(i,\beta)>(p,q)$. If $(\alpha,\beta)\in C_{i,j}(\sigma)$, then
denote by $c\sigma_{(i,j)}^{\alpha,\beta}\in S_n^2$ the involution
such that
\begin{equation*}
\Supp{c\sigma_{(i,j)}^{\alpha,\beta}}=
(\Supp{\sigma}\setminus\{(i,j)\})\cup\{(i,\beta),(\alpha,j)\}.
\end{equation*}

\exam{If $n=8$, $\sigma=(4,1)(8,2)(7,6)$, then $(3,5)\in
C_{8,2}(\sigma)$, $c\sigma_{(8,2)}^{3,5}=(4,1)(3,2)(8,5)(7,6)$, so
one has
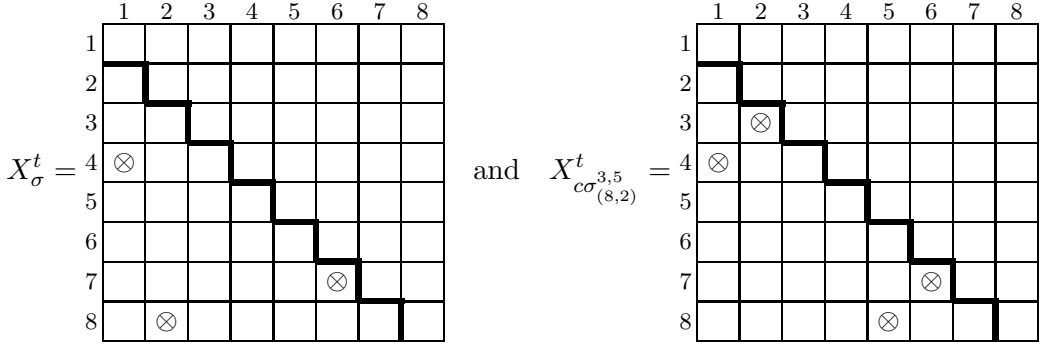
\begin{equation*}X_{\sigma}^t=
\mymatrix{ \pho& \pho& \pho& \pho& \pho& \pho& \pho& \pho\\
\Top{2pt}\Rt{2pt}\pho& \pho& \pho& \pho& \pho& \pho& \pho& \pho\\
\pho& \Top{2pt}\Rt{2pt}\pho& \pho& \pho& \pho& \pho& \pho& \pho\\
\otimes& \pho& \Top{2pt}\Rt{2pt}\pho& \pho& \pho& \pho& \pho& \pho\\
\pho& \pho& \pho& \Top{2pt}\Rt{2pt}\pho& \pho& \pho& \pho& \pho\\
\pho& \pho& \pho& \pho& \Top{2pt}\Rt{2pt}\pho& \pho& \pho& \pho\\
\pho& \pho& \pho& \pho& \pho& \Top{2pt}\Rt{2pt}\otimes& \pho& \pho\\
\pho& \otimes& \pho& \pho& \pho& \pho& \Top{2pt}\Rt{2pt}\pho& \pho\\
}\quad\text{and}\quad X_{c\sigma_{(8,2)}^{3,5}}^t=
\mymatrix{ \pho& \pho& \pho& \pho& \pho& \pho& \pho& \pho\\
\Top{2pt}\Rt{2pt}\pho& \pho& \pho& \pho& \pho& \pho& \pho& \pho\\
\pho& \Top{2pt}\Rt{2pt}\otimes& \pho& \pho& \pho& \pho& \pho& \pho\\
\otimes& \pho& \Top{2pt}\Rt{2pt}\pho& \pho& \pho& \pho& \pho& \pho\\
\pho& \pho& \pho& \Top{2pt}\Rt{2pt}\pho& \pho& \pho& \pho& \pho\\
\pho& \pho& \pho& \pho& \Top{2pt}\Rt{2pt}\pho& \pho& \pho& \pho\\
\pho& \pho& \pho& \pho& \pho& \Top{2pt}\Rt{2pt}\otimes& \pho& \pho\\
\pho& \pho& \pho& \pho& \otimes& \pho& \Top{2pt}\Rt{2pt}\pho& \pho\\
}\end{equation*}}

\sst Let $\sigma\in S_n^2$. Put $\mathrm{Near}(\sigma)=
N^+(\sigma)\cup N^-(\sigma)\cup N^0(\sigma)$, where
\begin{equation*}
\begin{split}
N^-(\sigma)&=\{\sigma_{(i,j)}^-,
(i,j)\in M(\sigma)\},\\
N^+(\sigma)&=\{c\sigma_{(i,j)}^{\alpha,\beta},
(i,j)\in\Supp{\sigma},(\alpha,\beta)\in C_{i,j}(\sigma)\},\\
N^0(\sigma)&=\{a\sigma_{(i,j)}^{(\alpha,\beta)},
(i,j)\in\Supp{\sigma},(\alpha,\beta)\in A_{i,j}(\sigma)\}\\
&\cup\{b\sigma_{(i,j)}^{(\alpha,\beta)},
(i,j)\in\Supp{\sigma},(\alpha,\beta)\in B_{i,j}(\sigma)\}\\
&\cup\{\sigma_{(i,j)}^{\to},
(i,j)\in\Supp{\sigma}\}\cup\{\sigma_{(i,j)}^{\uparrow},
(i,j)\in\Supp{\sigma}\}.\\
\end{split}
\end{equation*} \propp{Let $\sigma\in S_n^2$.
If $\tau\in\mathrm{Near}(\sigma)$, then\footnote{Cf. \cite[Lemma
3.16]{Melnikov2}.}
$\Omega_{\tau}\subset\overline{\Omega}_{\sigma}$.\label{prop:complex_closures}}{Denote
by $\overline{Z}^{\Cp}\subset\overline{Z}$ the closure of a subset
$Z\subseteq\nt^*$ in the complex topology. It's well-known that $\overline{\Omega}_{\sigma}^{\Cp}=\overline{\Omega}_{\sigma}$, so it's enough to prove that
$X_{\tau}^t\in\overline{\Omega}_{\sigma}^{\Cp}$.

i) First, assume $\tau=\sigma_{(i,j)}^-\in N^-(\sigma)$ for some
$(i,j)\in M(\sigma)$. Pick $\epsi\in\Cp^{\times}$ and put
\begin{equation*}
y_{\epsi}=x_{i,i}(\epsi).X_{\sigma}^t\in\Omega_{\sigma}.
\end{equation*}
Then
\begin{equation*}
(y_{\epsi})_{r,s}=\begin{cases}\epsi,&\text{if }r=i\text{
and }s=j,\\
(X_{\sigma}^t)_{r,s}=(X_{\tau}^t)_{r,s}&\text{otherwise}.
\end{cases}
\end{equation*}
Clearly, $y_{\epsi}\to X_{\tau}^t$ as $\epsi\to0$.

ii) Second, suppose $\tau=c\sigma_{(i,j)}^{\alpha,\beta}\in
N^+(\sigma)$ for some $(i,j)\in\Supp{\sigma}$, $(\alpha,\beta)\in
C_{i,j}(\sigma)$. Pick $\epsi\in\Cp^{\times}$ and put
\begin{equation*}
y_{\epsi}=x_{\alpha,i}(\epsi^{-1}).x_{j,\beta}(-\epsi^{-1}).x_{i,i}
(\epsi).X_{\sigma}^t\in\Omega_{\sigma}.
\end{equation*}
Then
\begin{equation*}
(y_{\epsi})_{r,s}=\begin{cases}1,&\text{if either }r=\alpha\text{
and }s=j,\text{ or }r=i\text{ and }s=\beta,\\
\epsi,&\text{if }r=i\text{
and }s=j,\\
(X_{\sigma}^t)_{r,s}=(X_{\tau}^t)_{r,s}&\text{otherwise}.
\end{cases}
\end{equation*}
Hence $y_{\epsi}\to X_{\tau}^t$ as $\epsi\to0$.

iii) Third, let $\tau=a\sigma_{(i,j)}^{(\alpha,\beta)}\in
N^0(\sigma)$ for some $(i,j)\in\Supp{\sigma}$, $(\alpha,\beta)\in
A_{i,j}(\sigma)$. Pick $\epsi\in\Cp^{\times}$ and put
\begin{equation*}
y_{\epsi}=x_{\beta,i}(\epsi^{-1}).x_{i,i}(\epsi).x_{\alpha,\alpha}
(-\epsi).X_{\sigma}^t\in\Omega_{\sigma}.
\end{equation*}
Then
\begin{equation*}
(y_{\epsi})_{r,s}=\begin{cases}1,&\text{if either }r=\beta\text{
and }s=j,\text{ or }r=\alpha\text{ and }s=i,\\
\epsi\text{ (resp. $-\epsi$)},&\text{if }r=i\text{
and }s=j\text{ (resp. }r=\alpha\text{ and $s=\beta$)},\\
(X_{\sigma}^t)_{r,s}=(X_{\tau}^t)_{r,s}&\text{otherwise}.
\end{cases}
\end{equation*}
Thus, $y_{\epsi}\to X_{\tau}^t$ as $\epsi\to0$.

iv) Now, if $\tau=b\sigma_{(i,j)}^{(\alpha,\beta)}\in N^0(\sigma)$
for some $(i,j)\in\Supp{\sigma}$, $(\alpha,\beta)\in
B_{i,j}(\sigma)$, then pick $\epsi\in\Cp^{\times}$, $\epsi\neq1$ and put
\begin{equation*}
y_{\epsi}=x_{j,\beta}(-\epsi^{-1}).x_{i,\alpha}(\epsi^{-1}).x_{\alpha,\alpha}
(\epsi).x_{i,i}(\epsi-\epsi^{-1}).X_{\sigma}^t\in\Omega_{\sigma}.
\end{equation*}
Then
\begin{equation*}
(y_{\epsi})_{r,s}=\begin{cases}1,&\text{if either }r=i\text{
and }s=\beta,\text{ or }r=\alpha\text{ and }s=j,\\
\epsi,&\text{if either }r=i\text{
and }s=j,\text{ or }r=\alpha\text{ and }s=\beta,\\
(X_{\sigma}^t)_{r,s}=(X_{\tau}^t)_{r,s}&\text{otherwise}.
\end{cases}
\end{equation*}
We see that $y_{\epsi}\to X_{\tau}^t$ as $\epsi\to0$.

v) Finally, suppose $\tau=\sigma_{(i,j)}^{\to}\in N^0(\sigma)$
(resp. $\tau=\sigma_{(i,j)}^{\uparrow}\in N^0(\sigma)$) for some
$(i,j)\in\Supp{\sigma}$. Pick $\epsi\in\Cp^{\times}$ and put
\begin{equation*}
y_{\epsi}=x_{j,m}(-\epsi^{-1}).x_{i,i}(\epsi).X_{\sigma}^t\in\Omega_{\sigma}
\text{ (resp. }
y_{\epsi}=x_{m,i}(\epsi^{-1}).x_{i,i}(\epsi).X_{\sigma}^t\in\Omega_{\sigma}).
\end{equation*}
Then
\begin{equation*}
(y_{\epsi})_{r,s}=\begin{cases}1,&\text{if }r=i\text{
and }s=m\text{ (resp. }r=m\text{ and }s=j\text{)},\\
\epsi,&\text{if }r=i\text{
and }s=j,\\
(X_{\sigma}^t)_{r,s}=(X_{\tau}^t)_{r,s}&\text{otherwise}.
\end{cases}
\end{equation*}
We conclude that $y_{\epsi}\to X_{\tau}^t$ as $\epsi\to0$. The
result follows.}

\medskip Things now are ready to prove the ``if'' direction of Theorem
\ref{mtheo} (see the next Section for the proofs of some technical but crucial
Lemmas). \propp{Let $\sigma$, $\tau$ be involutions in $S_n$. If
$\tau\leq^*\sigma$, then
$\Omega_{\tau}\subseteq\overline{\Omega}_{\sigma}$.\label{prop:only_if}}{Denote
$s(\sigma)=|\Supp{\sigma}|$ and put $L(\sigma)=L^+(\sigma)\cup
L^-(\sigma)\cup L^0(\sigma)$, where\footnote{Cf. \cite[Subsection
3.7]{Melnikov2}.}
\begin{equation*}\predisplaypenalty=0
\begin{split}
L^-(\sigma)&=\{\sigma'\in S_n^2\mid \sigma'\leq^*\sigma,\text{
}s(\sigma')<s(\sigma),\text{ and
if }\sigma'\leq^* w<^*\sigma,\text{ }s(w)<s(\sigma),\text{ then }w=\sigma'\},\\
L^+(\sigma)&=\{\sigma'\in S_n^2\mid \sigma'\leq^*\sigma,\text{
}s(\sigma')>s(\sigma),\text{ and
if }\sigma'\leq^* w<^*\sigma,\text{ then }w=\sigma'\},\\
L^0(\sigma)&=\{\sigma'\in S_n^2\mid \sigma'\leq^*\sigma,\text{
}s(\sigma')=s(\sigma),\text{ and if }\sigma'\leq^* w<^*\sigma,\text{
then }w=\sigma'\}.
\end{split}
\end{equation*}
Evidently, there exist involutions
$\sigma=w_1\geq^*w_2\geq^*\ldots\geq^*w_k=\tau$ such that
$w_{i+1}\in L(w_i)$ for all $1\leq i<k$, so we can assume $\tau\in
L(\sigma)$. But Lemmas \ref{lemm:hard_minus}, \ref{lemm:hard_ravno} and \ref{lemm:hard_plus}
show that $N^-(\sigma)=L^-(\sigma)$, $N^0(\sigma)=L^0(\sigma)$ and
$N^+(\sigma)=L^+(\sigma)$ respectively, so
$L(\sigma)=\mathrm{Near}(\sigma)$. Applying Proposition
\ref{prop:complex_closures}, we conclude the proof.}

\sst\label{sst:Bruhat} In this Subsection, we prove Theorem \ref{mcoro} (using technical Lemmas proved in the next Section). Let $\sigma$,~$\tau\in S_n^2$. Recall that
\begin{equation*}
\begin{split}
\tau\leq_B\sigma&\text{ if and only if }R(\dot\tau)\leq R(\dot\sigma),\\
\tau\leq^*\sigma&\text{ if and only if }R_{\tau}^*=R(\dot\tau)_{\mathrm{low}}\leq R_{\sigma}^*=R(\dot\sigma)_{\mathrm{low}},\\
\end{split}
\end{equation*}
so $\tau\leq_B\sigma$ implies $\tau\leq^*\sigma$ (see Subsection \ref{sst:Bruhat_order}).

In order to check that the converse holds, denote
\begin{equation}
\begin{split}
&L_B(\sigma)=\{\sigma'\in S_n^2\mid s'\leq_B\sigma\text{ and if }\sigma'\leq_Bw<_B\sigma,\text{  then }w=\sigma'\},\\
&L_*(\sigma)=\{\sigma'\in S_n^2\mid s'\leq^*\sigma\text{ and if }\sigma'\leq^*w<^*\sigma,\text{  then }w=\sigma'\}.
\end{split}
\end{equation}
Clearly, $L_*(\sigma)=L'(\sigma)\cup L^+(\sigma)\cup L^0(\sigma)$, where
\begin{equation*}
L'(\sigma)=\{\sigma'\in S_n^2\mid \sigma'\leq^*\sigma,\text{
}s(\sigma')<s(\sigma),\text{ and if }\sigma'\leq^* w<^*\sigma,\text{ then }w=\sigma'\}\subseteq L^-(\sigma).
\end{equation*}
(In general, $L'(\sigma)\subsetneq L^-(\sigma)$.) Put also
\begin{equation*}
N'(\sigma)=\{\sigma_{(i,j)}^-,
(i,j)\in M(\sigma)\text{ and }\sigma(m)\neq m\text{ for all }j\leq m\leq i\}\subseteq N^-(\sigma).
\end{equation*}

It follows from \cite[Theorem 5.2]{Incitti2} that $L_B(\sigma)=\mathrm{Near}'(\sigma)=N'({\sigma})\cup N^+(\sigma)\cup N^0(\sigma)$. But Lemmas~\ref{lemm:hard_ravno},~\ref{lemm:hard_plus}~and~\ref{lemm:hard_bis} show that $N^0(\sigma)=L^0(\sigma)$, $N^+(\sigma)=L^+(\sigma)$ and $N'(\sigma)=L'(\sigma)$ respectively, so $\mathrm{Near}'(\sigma)=L_*(\sigma)$. Hence the conditions $\sigma\geq^*\tau$ and $\sigma\geq_B\tau$ are equivalent; this proves Theorem~\ref{mcoro} and so concludes the proof of Theorem \ref{mtheo_0}. Furthermore, this gives the following combinatorial {de\-scrip\-tion} of Bruhat--Chevalley order on $S_n^2$:
\begin{equation*}
\tau\leq_B\sigma\text{ if and only if }R_{\tau}^*\leq R_{\sigma}^*.
\end{equation*}

\sect{Proofs of technical Lemmas\label{sect:proofs_Lemmas}}

\sst It turns out that the equalities
$L(\sigma)=\mathrm{Near}(\sigma)$ and $L_*(\sigma)=\mathrm{Near}'(\sigma)$ play a key role in the proofs of
Proposition~\ref{prop:only_if} and Theorem \ref{mcoro} respectively. The proofs of these equalities are
completely straightforward, but rather long. First, we will prove
that $\mathrm{Near}(\sigma)\subseteq L(\sigma)$. Obviously,
\begin{equation}
s(\tau)=s(\sigma)\pm1\text{ for all }\tau\in N^{\pm}(\sigma),\text{
and }s(\tau)=s(\sigma)\text{ for all }\tau\in
N^0(\sigma).\label{formula:nuzhnaya_dlina}
\end{equation}

Note that if $\sigma\geq^*w\geq^*\tau$, then
\begin{equation}
Y\cap\Supp{\sigma}=Y\cap\Supp{\tau}\Longrightarrow
Y\cap\Supp{w}=Y\cap\Supp{\sigma}=Y\cap\Supp{\tau}\label{formula:ravno_mezhdu}
\end{equation}
for all $Y\subseteq\Phi$ such that if $(a,b)\in Y$, $(c,d)\in\Phi$ and $(c,d)>(a,b)$, then $(c,d)\in Y$. Note also that
\begin{equation}
\sigma=\tau\Longleftrightarrow
R_{\sigma}^*=R_{\tau}^*,\label{formula:w1_w2_iff_R1_R2}
\end{equation}
and, moreover,
\begin{equation}
Y\cap\Supp{\sigma}=Y\cap\Supp{\tau}\Longleftrightarrow
(R_{\sigma}^*)_{r,s}=(R_{\tau}^*)_{r,s}\text{ for all }(r,s)\in
Y.\label{formula:R_cap_Y}
\end{equation}

Let $1\leq i,j\leq n$. It's very convenient to put
\begin{equation*}
\Ro_i=\{(i,s)\in\Phi\mid1\leq s<i\},\text{ }\Co_j=\{(r,j)\in\Phi\mid {j<r\leq n}\}.
\end{equation*}

\defi{The sets $\Ro_i$, $\Co_j$ are called the $i$th \emph{row} and the
$j$th \emph{column} of $\Phi$ respectively. Note that if $\sigma\in
S_n^2$, then
\begin{equation}
|\Supp{\sigma}\cap(\Ro_i\cup\Co_i)|\leq1\text{ for all }1\leq i\leq
n.\label{formula:ortog_supp}
\end{equation}}

\medskip\lemmp{Let $\sigma\in S_n^2$. One has\footnote{Cf.
\cite[Lemma 3.8]{Melnikov2}.} $N^-(\sigma)\subseteq
L^-(\sigma)$.\label{lemm:N_minus_in_L_minus}}{Suppose $\tau=\sigma_{(i,j)}^-$ for some $(i,j)\in
M(\sigma)$. By (\ref{formula:nuzhnaya_dlina}),
$s(\tau)=s(\sigma)-1<s(\sigma)$. Put
\begin{equation*}
Y=\Phi\setminus\{(p,q)\in\Phi\mid(p,q)\nleq(i,j)\}.
\end{equation*}
Then $\Supp{\tau}=\Supp{\tau}\cap Y=\Supp{\sigma}\cap Y$.

Now, assume that there exists $w\in S_n^2$ such that
$\tau\leq^*w<^*\sigma$ and $s(w)<s(\sigma)$. Then, by
(\ref{formula:ravno_mezhdu}),
$Y\cap\Supp{\sigma}=Y\cap\Supp{\tau}=Y\cap\Supp{w}$, so $s(w)\geq
s(\tau)=s(\sigma)-1$. Thus, $s(w)=s(\tau)$ and
$\Supp{w}=\Supp{\tau}$, so $w=\tau$.}

\medskip\lemmp{Let $\sigma\in S_n^2$. One
has\footnote{Cf. \cite[Lemmas 3.11--3.14]{Melnikov2}.}
$N^0(\sigma)\subseteq L^0(\sigma)$.\label{lemm:N_ravno_in_L_ravno}}{i) Let $\tau\in N^0(\sigma)$,
then, by (\ref{formula:nuzhnaya_dlina}), $s(\tau)=s(\sigma)$. First,
suppose $\tau=\sigma_{(i,j)}^{\uparrow}$ for some
$(i,j)\in\Supp{\sigma}$ (the case $\tau=\sigma_{(i,j)}^{\to}$ is
completely similar). Let
$\Supp{\tau}\setminus\Supp{\sigma}=\{(m,j)\}$. Put $Y=Y_0\cup Y_1$
and $\wt Y=\Phi\setminus Y$, where
\begin{equation*}
\begin{split}
&Y_0=\{(p,q)\in\Phi\mid(p,q)\nleq(i,j)\},\\
&Y_1=\{(p,q)\in\Phi\mid(p,q)\leq(m,j)\}.
\end{split}
\end{equation*}
Then $(R_{\sigma}^*)_{r,s}=(R_{\tau}^*)_{r,s}$ for all $(r,s)\in
Y$.

For example, let $n=8$, $i=7$, $j=2$, $m=5$. On the picture below
boxes from $Y_0$ are filled by~$0$'s, boxes from $Y_1$ are filled by
$1$'s, and boxes from $\wt Y$ are grey.
\begin{equation*}
\mymatrix{\pho& \pho& \pho& \pho& \pho& \pho& \pho& \pho\\
\Top{2pt}\Rt{2pt}0& \pho& \pho& \pho& \pho& \pho& \pho& \pho\\
\Rt{2pt}0& \Top{2pt}\Rt{2pt}1& \pho& \pho& \pho& \pho& \pho& \pho\\
\Rt{2pt}0& 1& \Top{2pt}\Rt{2pt}1& \pho& \pho& \pho& \pho& \pho\\
\Rt{2pt}0& 1& 1& \Top{2pt}\Rt{2pt}1& \pho& \pho& \pho& \pho\\
\Rt{2pt}0& \Top{2pt}\gray\pho& \Top{2pt}\gray\pho& \Top{2pt}\gray\pho& \Top{2pt}\Rt{2pt}\gray\pho& \pho& \pho& \pho\\
\Rt{2pt}0& \gray\pho& \gray\pho& \gray\pho& \gray\pho& \Top{2pt}\Rt{2pt}\gray\pho& \pho& \pho\\
0& \Top{2pt}0& \Top{2pt}0& \Top{2pt}0& \Top{2pt}0& \Top{2pt}0& \Top{2pt}\Rt{2pt}0& \pho\\
}\end{equation*}

\medskip Now, assume there exists $w\in S_n^2$ such that $\tau\leq^*
w<^*\sigma$. By (\ref{formula:w1_w2_iff_R1_R2}), it's enough to show
that $(R_{w}^*)_{r,s}=(R_{\tau}^*)_{r,s}$ for all $(r,s)\in\wt Y$,
because $(R_{\sigma}^*)_{r,s}=(R_{\tau}^*)_{r,s}=(R_{w}^*)_{r,s}$
for all $(r,s)\in Y$. Note that
$(R_{\tau}^*)_{r,s}=(R_{\sigma}^*)_{r,s}-1$ for all $(r,s)\in\wt Y$.
Further, by definition of $\sigma_{(i,j)}^{\uparrow}$ (see
Subsection \ref{sst:Near}), $\Supp{\sigma}\cap\wt Y=\{(i,j)\}$ and
$\Supp{\tau}\cap\wt Y=\emptyset$.

Since $\tau\leq^*w<^*\sigma$, there exists $(k,j)\in\Supp{w}$ such
that $m\leq k\leq i$. We claim that $\Supp{w}\cap(\wt
Y\setminus\Ro_i)=\emptyset$. Indeed, assume the converse holds,
i.e., there exists $(p,q)\in\Supp{w}\cap(\wt Y\setminus\Ro_i)$. Then
$m<p<i$. By definition of $\sigma_{(i,j)}^{\uparrow}$, there are no
$r$ such that $m<r<i$ and $\sigma(r)=r$. Hence $\Supp{\sigma}\cap
Y_0\cap(\Ro_p\cup\Co_p)\neq\emptyset$. By
(\ref{formula:ravno_mezhdu}), $\Supp{\sigma}\cap Y_0=\Supp{\tau}\cap
Y_0=\Supp{w}\cap Y_0$, so $\Supp{w}\cap
Y_0\cap(\Ro_p\cup\Co_p)\neq\emptyset$. But $(p,q)\in\wt Y$, hence
$|\Supp{\sigma}\cap (\Ro_p\cup\Co_p)|\geq2$. This contradicts
(\ref{formula:ortog_supp}). Thus, $\Supp{w}\cap(\wt
Y\setminus\Ro_i)=\emptyset$, as required.

In particular, either $k=i$ or $k=m$. If $k=i$, then
$\Supp{w}\cap\wt Y=\Supp{\sigma}\cap\wt Y=\{(i,j)\}$, so
$\Supp{w}\cap(\wt Y\cup Y_0)=\Supp{\sigma}\cap(\wt Y\cup Y_0)$, and,
by (\ref{formula:R_cap_Y}), $(R_{\sigma}^*)_{r,s}=(R_{w}^*)_{r,s}$
for all $(r,s)\in\wt Y\cup Y_0$. But
$(R_{\sigma}^*)_{r,s}=(R_{w}^*)_{r,s}$ for all $(r,s)\in Y_1$, hence
$R_{\sigma}^*=R_w^*$, and, by (\ref{formula:w1_w2_iff_R1_R2}),
$\sigma=w$, a contradiction. Thus, $k=m$. This means that
$\Supp{w}\cap\wt Y=\Supp{\tau}\cap\wt Y=\emptyset$, so
$\Supp{w}\cap(\wt Y\cup Y_0)=\Supp{\tau}\cap(\wt Y\cup Y_0)$, and,
by~(\ref{formula:R_cap_Y}), $(R_{\tau}^*)_{r,s}=(R_{w}^*)_{r,s}$ for
all $(r,s)\in\wt Y\cup Y_0$. But
$(R_{\tau}^*)_{r,s}=(R_{w}^*)_{r,s}$ for all $(r,s)\in Y_1$, hence
$R_{\tau}^*=R_w^*$, and, by (\ref{formula:w1_w2_iff_R1_R2}),
$\tau=w$, as required.

\medskip ii) Second, assume $\tau=b\sigma_{(i,j)}^{(\alpha,\beta)}$
for some $(i,j)\in\Supp{\sigma}$, $(\alpha,\beta)\in
B_{i,j}(\sigma)$. Then $(R_{\sigma}^*)_{r,s}=(R_{\tau}^*)_{r,s}$ for all $(r,s)\in Y$. Here we put $Y=Y_0\cup Y_1$ and $\wt Y=\Phi\setminus Y$,
where
\begin{equation*}
\begin{split}
&Y_0=\{(p,q)\in\Phi\mid(p,q)\nleq(\alpha,\beta)\},\\
&Y_1=\{(p,q)\in\Phi\mid(p,q)\leq(i,\beta)\text{ or
}(p,q)\leq(\alpha,j)\}.
\end{split}
\end{equation*}

\medskip For example, let $n=8$, $i=5$, $j=4$, $\alpha=7$, $\beta=1$.
On the picture below boxes from $Y_0$ are filled by~$0$'s, boxes
from $Y_1$ are filled by $1$'s, and boxes from $\wt Y$ are grey.
\begin{equation*}
\mymatrix{\pho& \pho& \pho& \pho& \pho& \pho& \pho& \pho\\
\Top{2pt}\Rt{2pt}1& \pho& \pho& \pho& \pho& \pho& \pho& \pho\\
1& \Top{2pt}\Rt{2pt}1& \pho& \pho& \pho& \pho& \pho& \pho\\
1& 1& \Top{2pt}\Rt{2pt}1& \pho& \pho& \pho& \pho& \pho\\
1& 1& 1& \Top{2pt}\Rt{2pt}1& \pho& \pho& \pho& \pho\\
\Top{2pt}\gray\pho& \gray\Top{2pt}\pho& \gray\Top{2pt}\Rt{2pt}\pho& 1& \Top{2pt}\Rt{2pt}1& \pho& \pho& \pho\\
\gray\pho& \gray\pho& \Rt{2pt}\gray\pho& 1& 1& \Top{2pt}\Rt{2pt}1& \pho& \pho\\
\Top{2pt}0& \Top{2pt}0& \Top{2pt}0& \Top{2pt}0& \Top{2pt}0& \Top{2pt}0& \Top{2pt}\Rt{2pt}0& \pho\\
}\end{equation*}

\medskip Now, assume there exists $w\in S_n^2$ such that $\tau\leq^*
w<^*\sigma$. By (\ref{formula:w1_w2_iff_R1_R2}), it's enough to show
that $(R_{w}^*)_{r,s}=(R_{\tau}^*)_{r,s}$ for all $(r,s)\in\wt Y$,
because $(R_{\sigma}^*)_{r,s}=(R_{\tau}^*)_{r,s}=(R_{w}^*)_{r,s}$
for all $(r,s)\in Y$. Note that
$(R_{\tau}^*)_{r,s}=(R_{\sigma}^*)_{r,s}-1$ for all $(r,s)\in\wt Y$.
Further, by definition of $b\sigma_{(i,j)}^{(\alpha,\beta)}$ (see
Subsection \ref{sst:Near}), $\Supp{\sigma}\cap\wt
Y=\{(\alpha,\beta)\}$ and $\Supp{\tau}\cap\wt Y=\emptyset$.

Since $\tau\leq^*w<^*\sigma$, there exists
$(k,\beta)\in\Supp{\sigma}$ such that $i\leq k\leq\alpha$. If
$k=\alpha$, then $w<^*\sigma$ follows $\Supp{w}\cap\wt
Y=\Supp{\sigma}\cap\wt Y=\{(\alpha,\beta)\}$, so $\Supp{w}\cap(\wt
Y\cup Y_0)=\Supp{\sigma}\cap(\wt Y\cup Y_0)$, and, by
(\ref{formula:R_cap_Y}), $(R_{\sigma}^*)_{r,s}=(R_{w}^*)_{r,s}$ for
all $(r,s)\in\wt Y\cup Y_0$. But
$(R_{\sigma}^*)_{r,s}=(R_{w}^*)_{r,s}$ for all $(r,s)\in Y_1$, hence
$R_{\sigma}^*=R_w^*$, and, by (\ref{formula:w1_w2_iff_R1_R2}),
$\sigma=w$, a contradiction. Thus, $k<\alpha$.

Similarly, there exists $(\alpha,l)\in\Supp{w}$ such that
$\beta<l\leq j$. If $k>i$ or $l<j$, then\linebreak
$(R_w^*)_{k,l}>(R_{\sigma}^*)_{k,l}$, which contradicts
$w<^*\sigma$. Hence $k=i$ and $l=j$, i.e., $(i,\beta)$ and
$(\alpha,j)$ belong to $\Supp{w}$. This implies $\Supp{w}\cap\wt
Y=\Supp{\tau}\cap\wt Y=\emptyset$, because if
$(p,q)\in\Supp{w}\cap\wt Y$, then
$(R_w^*)_{p,j}>(R_{\sigma}^*)_{p,j}$ and
$(R_w^*)_{i,q}>(R_{\sigma}^*)_{i,q}$, a contradiction. Thus,
$\Supp{w}\cap(\wt Y\cup Y_0)=\Supp{\tau}\cap(\wt Y\cup Y_0)$, and,
by (\ref{formula:R_cap_Y}), $(R_{\tau}^*)_{r,s}=(R_{w}^*)_{r,s}$ for
all $(r,s)\in\wt Y\cup Y_0$. But
$(R_{\tau}^*)_{r,s}=(R_{w}^*)_{r,s}$ for all $(r,s)\in Y_1$, hence
$R_{\tau}^*=R_w^*$, and, by (\ref{formula:w1_w2_iff_R1_R2}),
$\tau=w$, as required.

\medskip iii) Finally, suppose $\tau=a\sigma_{(i,j)}^{(\alpha,\beta)}$
for some $(i,j)\in\Supp{\sigma}$, $(\alpha,\beta)\in
A_{i,j}(\sigma)$. Then $(R_{\sigma}^*)_{r,s}=(R_{\tau}^*)_{r,s}$ for all $(r,s)\in Y$. Here we put $Y=Y_0\cup Y_1$ and $\wt Y=\Phi\setminus Y$,
where
\begin{equation*}
\begin{split}
&Y_0=\{(p,q)\in\Phi\mid(p,q)\nleq(i,j)\text{ and }
(p,q)\nleq(\alpha,\beta)\},\\
&Y_1=\{(p,q)\in\Phi\mid(p,q)\leq(\beta,j)\text{ or
}(p,q)\leq(\alpha,i)\}.
\end{split}
\end{equation*}

\medskip For example, let $n=8$, $i=6$, $j=2$, $\alpha=7$, $\beta=4$.
On the picture below boxes from $Y_0$ are filled by~$0$'s, boxes
from $Y_1$ are filled by $1$'s, and boxes from $\wt Y$ are grey.
\begin{equation*}
\mymatrix{\pho& \pho& \pho& \pho& \pho& \pho& \pho& \pho\\
\Top{2pt}\Rt{2pt}0& \pho& \pho& \pho& \pho& \pho& \pho& \pho\\
\Rt{2pt}0& \Top{2pt}\Rt{2pt}1& \pho& \pho& \pho& \pho& \pho& \pho\\
\Rt{2pt}0& 1& \Top{2pt}\Rt{2pt}1& \pho& \pho& \pho& \pho& \pho\\
\Rt{2pt}0& \gray\Top{2pt}\pho& \gray\Top{2pt}\pho& \gray\Top{2pt}\Rt{2pt}\pho& \pho& \pho& \pho& \pho\\
\Rt{2pt}0& \gray\pho& \gray\pho& \gray\pho& \gray\Top{2pt}\Rt{2pt}\pho& \pho& \pho& \pho\\
0& \Top{2pt}0& \Top{2pt}\Rt{2pt}0& \gray\pho& \gray\Rt{2pt}\pho& \Top{2pt}\Rt{2pt}1& \pho& \pho\\
0& 0& 0& \Top{2pt}0& \Top{2pt}0& \Top{2pt}0& \Top{2pt}\Rt{2pt}0& \pho\\
}\end{equation*}

\medskip Now, assume there exists $w\in S_n^2$ such that $\tau\leq^*
w<^*\sigma$. By (\ref{formula:w1_w2_iff_R1_R2}), it's enough to show
that $(R_{w}^*)_{r,s}=(R_{\tau}^*)_{r,s}$ for all $(r,s)\in\wt Y$,
because $(R_{\sigma}^*)_{r,s}=(R_{\tau}^*)_{r,s}=(R_{w}^*)_{r,s}$
for all $(r,s)\in Y$. Note that, by definition of
$a\sigma_{(i,j)}^{(\alpha,\beta)}$ (see Subsection \ref{sst:Near}),
$\Supp{\sigma}\cap\wt Y=\{(i,j),(\alpha,\beta)\}$ and
$\Supp{\tau}\cap\wt Y=\emptyset$.

Since $\tau\leq^*w<^*\sigma$, there exists $(k,j)\in\Supp{\sigma}$
such that $\beta\leq k\leq i$. We claim that\linebreak $\Supp{w}\cap
Y'=\emptyset$, where $Y'=\wt
Y\setminus\{(r,s)\in\Phi\mid(r,s)\geq(i,\beta)\})$. Indeed, assume
the converse holds, i.e., there exists $(p,q)\in\Supp{w}\cap Y'$.
Assume $\beta<p<i$ (the case $\beta<q<i$ is similar). By definition
of $a\sigma_{(i,j)}^{(\alpha,\beta)}$, there are no $r$ such that
$i<r<\beta$ and $\sigma(r)=r$. Hence $\Supp{\sigma}\cap
Y_0\cap(\Ro_p\cup\Co_p)\neq\emptyset$.
By~(\ref{formula:ravno_mezhdu}), $\Supp{\sigma}\cap
Y_0=\Supp{\tau}\cap Y_0=\Supp{w}\cap Y_0$, so $\Supp{w}\cap
Y_0\cap(\Ro_p\cup\Co_p)\neq\emptyset$. But $(p,q)\in\wt Y$, hence
$|\Supp{\sigma}\cap (\Ro_p\cup\Co_p)|\geq2$. This contradicts
(\ref{formula:ortog_supp}). Thus, $\Supp{w}\cap Y'=\emptyset$.

In particular, either $k=\beta$ or $k=i$, i.e., either
$(\beta,j)\in\Supp{w}$ or $(i,j)\in\Supp{w}$. Similarly, either
$(\alpha,\beta)\in\Supp{w}$ or $(\alpha,i)\in\Supp{w}$. If
$(i,j)\in\Supp{w}$, then, by (\ref{formula:ortog_supp}),
$(\alpha,i)\notin\Supp{w}$, so $(\alpha,\beta)\in\Supp{w}$.
Consequently, $\Supp{w}\cap\wt Y=\Supp{\sigma}\cap\wt
Y=\{(i,j),(\alpha,\beta)\}$, so $\Supp{w}\cap(\wt Y\cup
Y_0)=\Supp{\sigma}\cap(\wt Y\cup Y_0)$, and, by
(\ref{formula:R_cap_Y}), $(R_{\sigma}^*)_{r,s}=(R_{w}^*)_{r,s}$ for
all $(r,s)\in\wt Y\cup Y_0$. But
$(R_{\sigma}^*)_{r,s}=(R_{w}^*)_{r,s}$ for all $(r,s)\in Y_1$, hence
$R_{\sigma}^*=R_w^*$, and, by (\ref{formula:w1_w2_iff_R1_R2}),
$\sigma=w$, a contradiction. Hence $(\beta,j)\in\Supp{w}$.
By~(\ref{formula:ortog_supp}), $(\alpha,\beta)\notin\Supp{w}$, so
$(\alpha,i)\in\Supp{w}$. This implies $\Supp{w}\cap\wt
Y=\Supp{\tau}\cap\wt Y=\emptyset$. Thus, $\Supp{w}\cap(\wt Y\cup
Y_0)=\Supp{\tau}\cap(\wt Y\cup Y_0)$, and, by
(\ref{formula:R_cap_Y}), $(R_{\tau}^*)_{r,s}=(R_{w}^*)_{r,s}$ for
all $(r,s)\in\wt Y\cup Y_0$. But
$(R_{\tau}^*)_{r,s}=(R_{w}^*)_{r,s}$ for all $(r,s)\in Y_1$, hence
$R_{\tau}^*=R_w^*$, and, by (\ref{formula:w1_w2_iff_R1_R2}),
$\tau=w$, as required.}

\medskip\lemmp{Let $\sigma\in S_n^2$. One has\footnote{If
$\sigma\geq\tau$, then $s(\sigma)\geq s(\tau)$, so there are no
analogues to this Lemma in \cite{Melnikov2}.} $N^+(\sigma) \subseteq
L^+(\sigma)$.\label{lemm:N_plus_in_L_plus}}{Let $\tau\in
N^+(\sigma)$, then, by (\ref{formula:nuzhnaya_dlina}),
${s(\tau)=s(\sigma)+1>s(\sigma)}$. Suppose
$\tau=c\sigma_{(i,j)}^{\alpha,\beta}$ for some\linebreak
$(i,j)\in\Supp{\sigma}$, $(\alpha,\beta)\in C_{i,j}(\sigma)$. Then $(R_{\sigma}^*)_{r,s}=(R_{\tau}^*)_{r,s}$ for all $(r,s)\in Y$. Here we put
$Y=Y_0\cup Y_1$ and $\wt Y=\Phi\setminus Y$, where
\begin{equation*}
\begin{split}
&Y_0=\{(p,q)\in\Phi\mid(p,q)\nleq(i,j)\},\\
&Y_1=\{(p,q)\in\Phi\mid(p,q)\leq(\alpha,j)\text{ or
}(p,q)\leq(i,\beta)\}.
\end{split}
\end{equation*}

\medskip For example, let $n=8$, $i=7$, $j=2$, $\alpha=4$, $\beta=5$.
On the picture below boxes from $Y_0$ are filled by~$0$'s, boxes
from $Y_1$ are filled by $1$'s, and boxes from $\wt Y$ are grey.
\begin{equation*}
\mymatrix{\pho& \pho& \pho& \pho& \pho& \pho& \pho& \pho\\
\Top{2pt}\Rt{2pt}0& \pho& \pho& \pho& \pho& \pho& \pho& \pho\\
\Rt{2pt}0& \Top{2pt}\Rt{2pt}1& \pho& \pho& \pho& \pho& \pho& \pho\\
\Rt{2pt}0& 1& \Top{2pt}\Rt{2pt}1& \pho& \pho& \pho& \pho& \pho\\
\Rt{2pt}0& \gray\Top{2pt}\pho& \gray\Top{2pt}\pho& \gray\Top{2pt}\Rt{2pt}\pho& \pho& \pho& \pho& \pho\\
\Rt{2pt}0& \gray\pho& \gray\pho& \Rt{2pt}\gray\pho& \Top{2pt}\Rt{2pt}1& \pho& \pho& \pho\\
\Rt{2pt}0& \gray\pho& \gray\pho& \Rt{2pt}\gray\pho& 1& \Top{2pt}\Rt{2pt}1& \pho& \pho\\
0& \Top{2pt}0& \Top{2pt}0& \Top{2pt}0& \Top{2pt}0& \Top{2pt}0& \Top{2pt}\Rt{2pt}0& \pho\\
}\end{equation*}

\medskip Now, assume there exists $w\in S_n^2$ such that $\tau\leq^*
w<^*\sigma$. By (\ref{formula:w1_w2_iff_R1_R2}), it's enough to show
that $(R_{w}^*)_{r,s}=(R_{\tau}^*)_{r,s}$ for all $(r,s)\in\wt Y$,
because $(R_{\sigma}^*)_{r,s}=(R_{\tau}^*)_{r,s}=(R_{w}^*)_{r,s}$
for all $(r,s)\in Y$. Note that
$(R_{\tau}^*)_{r,s}=(R_{\sigma}^*)_{r,s}-1$ for all $(r,s)\in\wt Y$.
Further, by definition of $c\sigma_{(i,j)}^{\alpha,\beta}$ (see
Subsection \ref{sst:Near}), $\Supp{\sigma}\cap\wt Y=\{(i,j)\}$ and
$\Supp{\tau}\cap\wt Y=\emptyset$.

Since $\tau\leq^*w<^*\sigma$, there exists $(k,j)\in\Supp{\sigma}$
such that $\alpha\leq k\leq i$. If $k=i$, then $w<\sigma$ follows
$\Supp{w}\cap\wt Y=\Supp{\sigma}\cap\wt Y=\{(i,j)\}$, so
$$\Supp{w}\cap(\wt Y\cup Y_0)=\Supp{\sigma}\cap(\wt Y\cup Y_0)$$ and,
by (\ref{formula:R_cap_Y}), $(R_{\sigma}^*)_{r,s}=(R_{w}^*)_{r,s}$
for all $(r,s)\in\wt Y\cup Y_0$. But
$(R_{\sigma}^*)_{r,s}=(R_{w}^*)_{r,s}$ for all $(r,s)\in Y_1$, hence
$R_{\sigma}^*=R_w^*$, and, by (\ref{formula:w1_w2_iff_R1_R2}),
$\sigma=w$, a contradiction. Thus, $k<i$.

Similarly, there exists $(i,l)\in\Supp{w}$ such that $j<l<\beta$. We
claim that either $l\leq\alpha$ or $l=\beta$. Indeed, assume the
contrary, i.e., $\beta>l>\alpha$. By definition of
$c\sigma_{(i,j)}^{\alpha,\beta}$, there are no $s$ such that
$\alpha<s<\beta$ and $\sigma(r)=r$. Hence $\Supp{\sigma}\cap
Y_0\cap(\Ro_l\cup\Co_l)\neq\emptyset$. By
(\ref{formula:ravno_mezhdu}), $\Supp{\sigma}\cap Y_0=\Supp{\tau}\cap
Y_0=\Supp{w}\cap Y_0$, so $\Supp{w}\cap
Y_0\cap(\Ro_l\cup\Co_l)\neq\emptyset$. But $(i,l)\in\wt Y$, hence
$|\Supp{\sigma}\cap (\Ro_l\cup\Co_l)|\geq2$. This contradicts
(\ref{formula:ortog_supp}). Thus, either $l\leq\alpha$ or
$l=\beta$.

If $l=\alpha$, then $k=\alpha$ contradicts
(\ref{formula:ortog_supp}), so $k\neq\alpha$, i.e., $k>\alpha$. In this case, $(R_w^*)_{k,l}>(R_{\sigma}^*)_{k,l}$, which
contradicts $\sigma>^*w$. Hence either $l<\alpha$ or $l=\beta$. If
$l<\alpha$, then $(R_w^*)_{k,l}>(R_{\sigma}^*)_{k,l}$, as above.
This contradicts $\sigma>^*w$, so $l=\beta$. Similarly, $k=\alpha$,
Thus, $\Supp{w}\cap(\wt Y\cup Y_0)=\Supp{\tau}\cap(\wt Y\cup Y_0)$,
and, by~(\ref{formula:R_cap_Y}),
$(R_{\tau}^*)_{r,s}=(R_{w}^*)_{r,s}$ for all $(r,s)\in\wt Y\cup
Y_0$. But $(R_{\tau}^*)_{r,s}=(R_{w}^*)_{r,s}$ for all $(r,s)\in
Y_1$, hence $R_{\tau}^*=R_w^*$, and, by
(\ref{formula:w1_w2_iff_R1_R2}), $\tau=w$, as required.}

\sst In this Subsection, we will prove the most complicated parts of
Proposition \ref{prop:only_if}. Namely, we will show that
$L(\sigma)\subseteq\mathrm{Near}(\sigma)$ for all $\sigma\in S_n^2$, and so
these sets coincide. Note that
\begin{equation*}
w_0=(n,1)(n-1,2)\ldots(n-n_0+1,n_0),
\end{equation*}
where $n_0=[n/2]$, is the maximal element
of $S_n^2$ with respect to the partial order $\leq^*$ on $S_n^2$.

\medskip\lemmp{Let $\sigma\in S_n^2$. Then $L^-(\sigma)=N^-(\sigma)$.
\label{lemm:hard_minus}}{By
Lemma \ref{lemm:N_minus_in_L_minus}, it's enough to check that
$L^-(\sigma)\subseteq N^-(\sigma)$. We must show that
\begin{equation}
\begin{split}
&\text{if }
\tau\leq^*\sigma\text{ and }s(\tau)<s(\sigma),\\
&\text{then there exists }\sigma'\in N^-(\sigma)
\text{ such that }\tau\leq^*\sigma'<^*\sigma.
\end{split}\label{formula:proof_minus}
\end{equation}
We will proceed by induction on $n$
(for $n=1$, there is nothing to prove). The proof is rather long, so we split it into six
steps.

\medskip i) Let $\sigma=w_0$, the maximal element of $S_n^2$ with respect to
$\leq^*$, $s(\tau)<s(w_0)$ and $\tau<^* w_0$. Let $\sigma'$ be the involution
such that $\Supp{\sigma'}=\Supp{\sigma}\setminus\{(n_0,n-n_0+1)\}$, where
$n_0=[n/2]$. Since $(n_0,n-n_0+1)\in M(w_0)$, $\sigma'\in N^-(w_0)$. But
$s(\tau)<s(w_0)$ implies
\begin{equation*}
(R_{\tau}^*)_{i,j}\leq|\Supp{\tau}|\leq|\Supp{w_0}|-1\leq(R_{\sigma'}^*)_{i,j}
\end{equation*}
for all $(i,j)\in\Phi$, so $\sigma'\geq^*\tau$. Therefore, we may also use the second (downward) induction on $\leq^*$.

\medskip ii) Let $\sigma=(i_1,j_1)\ldots(i_s,j_s)\in S_n^2$, $\sigma<^*w_0$,
$\tau=(p_1,q_1)\ldots(p_t,q_t)<^*\sigma$ and $s>t$.
Consider the following conditions:
\begin{equation}
\begin{split}
&\text{a) There exists }k\leq n_0\text{ such that }
j_l=l\text{ for all }1\leq l\leq k\text{ and either }i_k=k+1
\text{ or }k=s.\\
&\text{b) There exists }d\leq k\text{ such that }q_l=l\text{ for all }
1\leq l\leq d\text{ and either }q_{d+1}>k
\text{ or }d=t.\\
&\text{c) }(i_l,j_l)=(i_l,l)>(i_{l+1},j_{l+1})=
(i_{l+1},l+1),\text{ i.e., }i_l>i_{l+1},\text{ for all }1\leq l\leq k-1.\\
&\text{d) }(p_l,q_l)=(p_l,l)>(p_{l+1},q_{l+1})=
(p_{l+1},l+1),\text{ i.e., }p_l>p_{l+1},\text{ for all }1\leq l\leq d-1.\\
&\text{e) }(i_l,j_l)=(i_l,l)\geq(p_l,q_l)=
(p_l,l),\text{ i.e., }i_l\geq p_l,\text{ for all }1\leq l\leq d.
\end{split}\label{formula:remained_case_minus}
\end{equation}
Pick $r\leq s$. Define $\sigma_r$ and $\tau_r$ by putting
$\Supp{\sigma_r}=\Supp{\sigma}\cap\wt\Co_r$,
$\Supp{\tau_r}=\Supp{\tau}\cap\wt\Co_r$, where $\wt\Co_r=
\bigcup_{l\leq r}\Co_l$. We claim that
\begin{equation}
\begin{split}
&\text{if (\ref{formula:proof_minus}) holds for all }
\sigma,\tau\text{ satisfying (\ref{formula:remained_case_minus})},\\
&\text{then (\ref{formula:proof_minus})
holds for all }\sigma,\tau\in S_n^2.
\end{split}\label{formula:remained_minus_r}
\end{equation}
Clearly,
it's enough to prove that if (\ref{formula:proof_minus}) holds for all
$\sigma $, $\tau$ satisfying (\ref{formula:remained_case_minus}), and $\sigma_r$, $\tau_r$
don't satisfy~(\ref{formula:remained_case_minus}) for some $1\leq r\leq s$,
then (\ref{formula:proof_minus}) holds for $\sigma$,
$\tau$. We will proceed by induction on $r$. Evidently,
there exist $w_1=\sigma$, $w_2$, $\ldots$, $w_h=\tau\in S_n^2$
such that $w_1>^*w_2>^*\ldots>^*w_h$ and $w_{i+1}\in L^-(w_i)$
for all $1\leq i<h$, so we may assume $\tau\in L^-(\sigma)$.

The base $r=1$ is clear. Indeed, if $j_1>i_1$, then it follows from $\tau<^*\sigma$
that $q_1>1$, so $\sigma$, $\tau$ belong to $\wt S_{n-1}=
\{w\in S_n\mid w(1)=1\}\cong S_{n-1}$, and
we can use the first inductive assumption. Hence $\sigma_1$
satisfies~(\ref{formula:remained_case_minus}a); $\tau_1$~satisfies
(\ref{formula:remained_case_minus}b) automatically. In the case $r=1$ conditions
(\ref{formula:remained_case_minus}c) and (\ref{formula:remained_case_minus}d) are
trivial, so it remains to check (\ref{formula:remained_case_minus}e). But if $q_1=1$
and $p_1>i_1$, then $\tau\nless^*\sigma$, a contradiction. Thus, $\sigma_1$,
$\tau_1$ satisfy~(\ref{formula:remained_case_minus}e), as~required.

\medskip iii) Now, suppose $1\leq r\leq s$ and $\sigma_r$, $\tau_r$ satisfy
(\ref{formula:remained_case_minus}). To perform the induction step,
we must prove that either $\sigma_{r+1}$, $\tau_{r+1}$ satisfy
(\ref{formula:remained_case_minus}), too, or (\ref{formula:proof_minus}) holds
for $\sigma$, $\tau$. This is trivially true if
$i_k=k+1$ for some $k\leq r$, so we may assume that $i_l>l+1$ for all
$1\leq l\leq r$. Since $r\leq s$ and $\sigma_r$ satisfies
(\ref{formula:remained_case_minus}a), $j_r=r$. First,
consider the case when $p_l>l+1$ for all $l\leq r_0=\min\{r,s(\tau_r)\}$.
Suppose $j_{r+1}>r+1$, i.e., $\Supp{\sigma}\cap\Co_{r+1}=\emptyset$. Put
\begin{equation*}
\begin{split}
&\wt\sigma=P_r(\sigma)=(i_1,2)(i_2,3)\ldots(i_r,r+1)(i_{r+1},j_{r+1})\ldots(i_s,j_s),\\
&\wt\tau=P_r(\tau)=(p_1,2)(p_2,3)\ldots(p_{r_0},r_0+1)(p_{r_0+1},q_{r_0+1})\ldots(p_t,q_t).
\end{split}
\end{equation*}
Here the map $P_r\colon S_n^2\to S_n$ is defined by the following rule: if
$\eta=(a_1,b_1)\ldots(a_m,b_m)\in S_n^2$, then
\begin{equation*}
P_r(\eta)=
(a_1,b_1+1)\ldots(a_{r},b_z+1)(a_{r+1},b_{z+1})\ldots(a_m,b_m),
\end{equation*}
where $b_z\leq r<b_{z+1}$.
Note that, in general, $P_r(\eta)$ is \emph{not} an involution.

Clearly, $\wt\sigma$ and
$\wt\tau$ are \emph{involutions} in $S_n$. Indeed, it
suffice to check that (\ref{formula:ortog_supp}) holds for
$\Supp{\wt\tau}$ (for $\Supp{\wt\sigma}$, there is nothing to check,
because $r+1<j_{r+1}$). But if
$r_0<r$, then $r_0+1<r+1\leq q_{r_0+1}$, because $\tau_r$ satisfies
(\ref{formula:remained_case_minus}b), so (\ref{formula:ortog_supp}) holds
for $\Supp{\wt\tau}$. On the other hand, if $r=r_0$, then $\tau<^*\sigma$
yields $\Supp{\tau}\cap\Co_{r+1}=\emptyset$ (if the converse holds, then
$(R_{\sigma}^*)_{r+2,r+1}=r<r+1=(R_{\tau}^*)_{r+2,r+1}$, a contradiction).
Hence $r_0+1<q_{r_0+1}$ and (\ref{formula:ortog_supp}) holds
for $\Supp{\wt\tau}$. Thus, $\wt\sigma,\wt\tau\in S_n^2$.

Further, they belong to
$\wt S_{n-1}$ and $s(\wt\tau)=t<s=s(\wt\sigma)$, so, by the first induction hypothesis,
there exists $\wt w\in\wt S_{n-1}$ such that $\wt\sigma>^*\wt w\geq^*\wt\tau$
and $\wt w=\wt\sigma_{(i,j)}^-\in N^-(\wt\sigma)$ for some $(i,j)\in M(\wt\sigma)$.
Since $\sigma_r$~satisfies (\ref{formula:remained_case_minus}a),
$M(\wt\sigma)\cap\Co_{l+1}=\emptyset$ for all $1\leq l\leq r-1$, so $j\geq r+1$.
If $j>r+1$, then $(i,j)\in M(\sigma)$ and $w=\sigma_{(i,j)}^-\geq^*\tau$.
On the other hand, if $j=r+1$, then $(i,r)\in M(\sigma)$ and
$w=\sigma_{i,r}^-\geq^*\tau$. In both cases, $w\in N^-(\sigma)$
and $w\geq^*\tau$, as required.

\medskip iv) Then, suppose $\sigma_r$, $\tau_r$ satisfy
(\ref{formula:remained_case_minus}), $i_l>l+1$ for all
$1\leq l\leq r$, $p_l>l+1$ for all $l\leq r_0=\min\{r,s(\tau_r)\}$, but
$j_{r+1}=r+1$, i.e., $\Supp{\sigma}\cap\Co_{r+1}\neq\emptyset$. If
$i_{r+1}>i_r$, then put
\begin{equation*}
\sigma_0=(i_1,1)\ldots(i_{r-1},r-1)
(i_{r+1},r)(i_r,r+1)(i_{r+2},j_{r+2})\ldots(i_s,j_s).
\end{equation*}
Then $\sigma_0>^*\sigma>^*\tau$ and $s(\sigma_0)=s>t=s(\tau)$, so,
by the second inductive assumption, there exists $w_1\in S_n$ such
that $w_1\geq^*\tau$ and $w_1=(\sigma_0)_{(i,j)}^-\in N^-(\sigma_0)$
for some $(i,j)\in M(\sigma_0)$. Since $\sigma_r$ satisfies
(\ref{formula:remained_case_minus}a),
$M(\sigma_0)\cap\Co_l=\emptyset$ for all $1\leq l\leq r-1$, so
$j\geq r+1$. If $j>r+1$, then $(i,j)\in M(\sigma)$ and
$w=\sigma_{(i,j)}^-\geq^*\tau$.

On the other hand, assume $j=r+1$. If $\Supp{\tau}\cap\Co_r=\emptyset$,
i.e., $s(\tau_r)<r$,
then $w=\sigma_{(i,r)}^-\geq^*\tau$. If $\Supp{\tau}\cap\Co_r\neq\emptyset$,
i.e., $q_r=r$, then set $\sigma_1$, $\tau_1$ to be the involutions
such that
\begin{equation*}
\Supp{\sigma_1}=\Supp{\sigma}\setminus\Co_r,\text{ }\Supp{\tau_1}=\Supp{\tau}\setminus\Co_r.
\end{equation*}
Put also
$\wt\sigma_1=P_{r-1}(\sigma_1)$ and $\wt\tau_1=P_{r-1}(\tau_1)$.
We see that $\wt\sigma_1,\wt\tau_1\in S_n^2$ and $\wt\sigma_1,\wt\tau_1\in\wt S_{n-1}$.
Moreover, $s(\wt\sigma_1)=s-1>t-1=s(\wt\tau_1)$ and $\wt\sigma_1>^*\wt\tau_1$ (if
$\wt\sigma_1=\wt\tau_1$, then $s=t$, a~contradiction). Hence, by the first
induction hypothesis, there exists $\wt w_1\in S_n^2$ such that $\wt w_1\geq^*\wt\tau_1$
and $\wt w_1=(\wt\sigma_1)_{(a,b)}^-\in N^-(\wt\sigma_1)$ for some
$(a,b)\in M(\wt\sigma_1)$. Since $j_r=r$, $q_r=r$ and $\sigma_r$, $\tau_r$ satisfy
(\ref{formula:remained_case_minus}), $M(\wt\sigma_1)\cap\Co_{l+1}=\emptyset$
for all $1\leq l\leq r-1$, so $b\geq r+1$. Hence $(a,b-1)\in M(\sigma)$
and $w=\sigma_{(a,b-1)}^-\geq^*\tau$. We conclude that if $\sigma_{r+1}$
doesn't satisfy~(\ref{formula:remained_case_minus}c), then (\ref{formula:proof_minus})
holds for $\sigma$, $\tau$.

Next, suppose $i_{r+1}<i_r$, but $r_0=r$, $q_{r+1}=r+1$ and $p_{r+1}>i_{r+1}$. In
this case, put
\begin{equation*}
\tau_0=(p_1,1)\ldots(p_{r-1},r-1)(p_{r+1},r)(p_r,r+1)(p_{r_2},q_{r+2})
\ldots(p_t,q_t).
\end{equation*}
Then $\sigma>^*\tau_0>^*\tau$, so $\tau\notin L^-(\sigma)$.
The cases $i_{r+1}<i_r$, $r_0=r$,
$q_{r+1}=r+1$, $p_r<p_{r+1}<i_{r+1}$, and $i_{r+1}<i_r$, $r_0<r$,
$\Supp{\tau}\cap\Co_{r+1}\neq\emptyset$ are similar.
Namely, if $i_{r+1}<i_r$, $r_0=r$,
$q_{r+1}=r+1$, $p_r<p_{r+1}<i_{r+1}$, then we define $\tau_0$
as above, and if $i_{r+1}<i_r$, $r_0<r$,
$\Supp{\tau}\cap\Co_{r+1}=\{(p_{r_0+1},r+1)\}\neq\emptyset$, then we put
\begin{equation*}
\tau_0=(p_1,1)\ldots(p_{r_0},r_0)(p_{r_0+1},r)(p_{r_0+2},q_{r_0+2})
\ldots(p_t,q_t).
\end{equation*}
In both cases, $\sigma>^*\tau_0>^*\tau$, so $\tau\notin L^-(\sigma)$.
We conclude that if
$\tau_{r+1}$ doesn't satisfy (\ref{formula:remained_case_minus}d), or
$\sigma_{r+1}$, $\tau_{r+1}$ don't
satisfy~(\ref{formula:remained_case_minus}e), then (\ref{formula:proof_minus})
holds for $\sigma$, $\tau$.

\medskip v) To prove (\ref{formula:remained_minus_r}), it remains to
consider the case when $\sigma_r$, $\tau_r$ satisfy
(\ref{formula:remained_case_minus}),
$i_l>l+1$ for all $1\leq l\leq r$, but $p_d=d+1$ for some $d\leq r$.
(And, consequently, $r_0=s(\tau_r)=d$.) Suppose $j_{r+1}>r+1$,
i.e., $\Supp{\sigma}\cap\Co_{r+1}=\emptyset$. Let $\tau_0$ be the
involution such that $\Supp{\tau_0}=\Supp{\tau}\setminus\{(d+1,d)\}$;
put also $\wt\sigma=P_r(\sigma)$, $\wt\tau_0=P_r(\tau_0)$. Then
$\wt\sigma$, $\wt\tau_0\in\wt S_{n-1}\cap S_n^2$, $\wt\sigma>^*\wt\tau_0$
and $s(\wt\sigma)=s>t>t-1=s(\wt\tau_0)$. Hence, by the first inductive
assumption, there exists $\wt w=\wt\sigma_{(i,j)}^-\in N^-(\wt\sigma)$ such that
$\wt w\geq^*\wt\tau_0$. Since $\sigma_r$ satisfies~(\ref{formula:remained_case_minus}a),
$M(\wt\sigma)\cap\Co_{l+1}=\emptyset$ for all $1\leq l\leq r-1$, so $j\geq r+1$. If
$j>r+1$, then $(i,j)\in M(\sigma)$ and $w=\sigma_{(i,j)}^-\geq^*\tau$. Similarly, if $j=r+1$ and $r>d$, then $(i_r,r)\in M(\sigma)$ and $w=\sigma_{(i_r, r)}^-\leq^*\tau$.

On the other hand, if $j=r+1$ and $r=d$, then put $\wt\sigma_1=P_d(\sigma_1)$, where $\sigma_1$ is the involution
such that $\Supp{\sigma_1}=\Supp{\sigma}\setminus\{(i_d,d)\}$. In this case, $\wt\sigma_1$,
$\wt\tau_0\in\wt S_{n-1}\cap S_n^2$, $s(\wt\sigma_1)=s-1>t-1=s(\wt\tau_0)$ and $\wt\sigma_1>^*\wt\tau_0$. Whence the first induction hypothesis
shows that there exists $\wt w_1=(\sigma_1)_{(\alpha,\beta)}^-\in N^-(\wt\sigma_1)$ such that
$\wt w_1\geq^*\wt\tau_0$. Since $\sigma_r$, $\tau_r$ satisfy (\ref{formula:remained_case_minus}), $\beta>r=d$ (and so $\beta>r+1$).
Thus, $(\alpha,\beta)\in M(\sigma)$ and $w=\sigma_{(\alpha,\beta)}^-\geq^*\tau$.

Now, suppose $j_{r+1}=r+1$, but $i_{r+1}>i_r$. Arguing as on the
previous step, we conclude that there exists $w\in N^-(\sigma)$ such
that $w\geq^*\tau$. Thus, if $\sigma_{r+1}$ doesn't satisfy~(\ref{formula:remained_case_minus}a), then (\ref{formula:proof_minus}) holds
for $\sigma$, $\tau$. Next, suppose $j_{r+1}=r+1$, $i_{r+1}<i_r$, but $\tau_{r+1}$ doesn't satisfy
(\ref{formula:remained_case_minus}c) or (\ref{formula:remained_case_minus}d), i.e.,
$\Supp{\tau}\cap\Co_{r+1}=\{(p,r+1)\}\neq\emptyset$ (and,
consequently, $r>d$, since (\ref{formula:ortog_supp}) holds for
$\tau$). If $r>d+1$, then put
\begin{equation*}
\tau_2=(p_1,1)\ldots(p_{d-1},d-1)(d+1,d)(p,r)\ldots(p_t,q_t),
\end{equation*}
i.e., $\Supp{\tau_2}=(\Supp{\tau}\setminus\{(p,r+1)\})
\cup\{(p,r)\}$. Then $\sigma>^*\tau_2>^*\tau$. Indeed,
the last inequality is evident, and the first one follows
from $\sigma>^*\tau$ and the fact that $\sigma_r$, $\tau_r$
satisfy (\ref{formula:remained_case_minus}). Thus,
$\tau\notin L^-(\sigma)$. Similarly, if $r=d+1$, then put
\begin{equation*}
\tau_3=(p_1,1)\ldots(p_{d-1},d-1)(p,d)(p_{d+2},q_{d+2})\ldots(p_t,q_t),
\end{equation*}
i.e., $\Supp{\tau_3}=(\Supp{\tau}\setminus\{(d+1,d),(p,r+1)\})
\cup\{(p,d)\}$. In this case, $\sigma>^*\tau_3>^*\tau$, so
$\tau\notin L^-(\sigma)$. The proof of (\ref{formula:remained_minus_r})
is complete.

\medskip vi) Now, we may assume without loss of generality that $\sigma$, $\tau$
satisfy (\ref{formula:remained_case_minus}). If $i_l>l+1$ for all
$1\leq l\leq k$, then, by (\ref{formula:remained_case_minus}a) and
(\ref{formula:remained_case_minus}b), $$k=s=s(\sigma)\geq d=t=s(\tau).$$
If $q_d>d+1$, then put $\wt\sigma=P_k(\sigma)$,
$\wt\tau=P_k(\tau)=P_d(\tau)$. As above, $\wt\sigma$, $\wt\tau\in
\wt S_{n-1}\cap S_n^2$, so, by the first induction hypothesis,
there exists $\wt w=\wt\sigma_{(i,j)}^-$ such that $\wt w\geq^*\wt\tau$.
Since $\sigma$, $\tau$
satisfy (\ref{formula:remained_case_minus}), $M(\wt\sigma)\cap\Co_{l+1}
=\emptyset$ for all $1\leq l\leq k$, so $(i,j)\in M(\sigma)$ and
$\tau\leq^*w=\sigma_{(i,j)}^-\in N^-(\sigma)$. On the other hand, if $q_d=d+1$
for some $d\leq k$, then, arguing as
on the previous step, we conclude that such an involution $w$ exists.

Finally, assume $i_k=k+1$ for some $k\leq s$. Clearly, $(k+1, k)\in
M(\sigma)$. If $\Supp{\tau}\cap\Co_k=\emptyset$, i.e., $d<k$, then
$\tau\leq^*w=\sigma_{(k+1,k)}^-\in N^-(\sigma)$. At the contrary, if
$d=k$, i.e., $(k+1,k)\in\Supp{\tau}$, then we derive an existence of
$w$ as on the previous step. The proof is complete.}

\medskip\lemmp{Let $\sigma\in S_n^2$. Then $L^0(\sigma)=N^0(\sigma)$.
\label{lemm:hard_ravno}}{By
Lemma \ref{lemm:N_ravno_in_L_ravno}, it's enough to check that
$L^0(\sigma)\subseteq N^0(\sigma)$. By definition,
\begin{equation*}
\begin{split}
&L^0(\sigma)=\{\sigma'\in S_n^2\mid \sigma'\leq^*\sigma,\text{
}s(\sigma')=s(\sigma),\text{ and if }\sigma'\leq^* w<^*\sigma,\text{
then }w=\sigma'\}\\
&\subseteq\wt L^0(\sigma)=\{\sigma'\in S_n^2\mid \sigma'\leq^*\sigma,\text{
}s(\sigma')=s(\sigma),\text{ and if }\sigma'<^*w<^*\sigma,
\text{ then }s(w)>s(\sigma)\},\\
\end{split}
\end{equation*}
so it suffice to show that $\wt L^0(\sigma)\subseteq N^0(\sigma)$, i.e.,
\begin{equation}\predisplaypenalty=0
\begin{split}
&\text{if }
\tau\leq^*\sigma\text{ and }s(\tau)=s(\sigma),\\
&\text{then there exists }\sigma'\in N^0(\sigma)\cup N^-(\sigma)
\text{ such that }\tau\leq^*\sigma'<^*\sigma.
\end{split}\label{formula:proof_ravno}
\end{equation}
We will proceed by induction on $n$
(for $n=1$, there is nothing to prove). The proof is rather long, so we split it into
seven steps.

\medskip i) Let $\sigma=w_0$, the maximal element of $S_n^2$ with respect to
$\leq^*$, $s(\tau)=s(w_0)$ and $\tau<w_0$. Let $n_0=[n/2]$ and
$h=\min\{s\mid1\leq s\leq n_0\text{ and }(n-s+1,s)\notin\Supp{\tau}\}$. (If
$h$ doesn't exist, then $R_{\sigma}^*=R_{\tau}^*$ and, by (\ref{formula:w1_w2_iff_R1_R2}),
$\tau=\sigma$, a contradiction.) If $s<n_0$, then $(R_{\tau}^*)_{n-s+1,s}\leq s-1$,
hence $\tau\leq^*w=b(w_0)_{(n-s,s+1)}^{(n-s+1,s)}\in N^0(w_0)$
(clearly, $(n-s+1,s)\in B_{n-s,s+1}(w_0)$), because
\begin{equation*}
(R_w^*)_{i,j}=\begin{cases}s-1,&\text{if }i=n-s+1\text{ and }j=s,\\
(R_{w_0}^*)_{i,j}&\text{otherwise}.
\end{cases}
\end{equation*}
On the other hand, if $s=n_0$, then $n$ is odd and either
\begin{equation*}
\begin{split}
&\Supp{\tau}=(\Supp{w_0}\setminus\{(n_0+2,n_0)\})
\cup\{(n_0+1,n_0)\}\text{ or}\\
&\Supp{\tau}=(\Supp{w_0}\setminus\{(n_0+2,n_0)\})\cup\{(n_0+2,n_0+1)\}.
\end{split}
\end{equation*}
If the first case occurs, then we put $w=(w_0)_{(n_0+2,n_0)}^{\uparrow}$,
and if the second case occurs, we put $w=(w_0)_{(n_0+2,n_0)}^{\to}$.
In both cases, $w\in N^0(w_0)$ and $w\geq^*\tau$ (in fact, $w=\tau$).
Therefore, we may also use the~second (downward) induction on $\leq^*$.

\medskip ii) Let $\sigma=(i_1,j_1)\ldots(i_t,j_t)\in S_n^2$, $\sigma<^*w_0$,
$\tau=(p_1,q_1)\ldots(p_t,q_t)<^*\sigma$.
Consider the following conditions (cf. (\ref{formula:remained_case_minus})):
\begin{equation}
\begin{split}
&\text{a) There exists }k\leq t\text{ such that }
i_l>l+1\text{ for all }1\leq l\leq k-1.\\
&\text{b) There exists }d\leq k\text{ such that }q_l=l\text{ for all }
1\leq l\leq d\\
&\hphantom{\text{b) }}\text{ and either }d<k,\text{ }q_{d+1}>k,
\text{ or }d=k,\text{ }p_d=d+1.\\
&\text{c) }j_l=l\text{ for all }1\leq l\leq d,\text{ and if }
d<k,\text{ then either }i_k=j_k+1\text{ or }k=t.\\
&\text{d) }(p_l,q_l)=(p_l,l)>(p_{l+1},q_{l+1})=
(p_{l+1},l+1),\text{ i.e., }p_l>p_{l+1},\text{ for all }1\leq l\leq d-1.\\
&\text{e) }(i_l,j_l)=(i_l,l)\geq(p_l,q_l)=
(p_l,l),\text{ i.e., }i_l\geq p_l,\text{ for all }1\leq l\leq d.
\end{split}\label{formula:remained_case_ravno}
\end{equation}
Pick $r\leq t$. Define $\sigma_r$ and $\tau_r$ as in the previous Lemma. We claim that
\begin{equation}
\begin{split}
&\text{if (\ref{formula:proof_ravno}) holds for all }
\sigma,\tau\text{ satisfying (\ref{formula:remained_case_ravno})},\\
&\text{then (\ref{formula:proof_ravno})
holds for all }\sigma,\tau\in S_n^2.
\end{split}\label{formula:remained_ravno_r}
\end{equation}
Clearly,
it's enough to prove that if (\ref{formula:proof_ravno}) holds for all
$\sigma $, $\tau$ satisfying (\ref{formula:remained_case_ravno}), and $\sigma_r$, $\tau_r$
don't satisfy~(\ref{formula:remained_case_ravno}) for some $1\leq r\leq t$,
then (\ref{formula:proof_ravno}) holds for $\sigma$,
$\tau$.

We will proceed by induction on $r$. Evidently,
there exist $w_1=\sigma$, $w_2$, $\ldots$, $w_z=\tau\in S_n^2$
such that $w_1>^*w_2>^*\ldots>^*w_z$ and $w_{i+1}\in\wt L^0(w_i)$
for all $1\leq i<z$, so we may assume $\tau\in\wt L^0(\sigma)$. The base $r=1$
is evident (see step ii) of the proof of Lemma \ref{lemm:hard_minus}).

\medskip iii) Suppose $1\leq r\leq t$ and $\sigma_r$, $\tau_r$ satisfy
(\ref{formula:remained_case_ravno}). To perform the induction step,
we must prove that either $\sigma_{r+1}$, $\tau_{r+1}$ satisfy
(\ref{formula:remained_case_ravno}), too, or
(\ref{formula:proof_ravno}) holds for $\sigma$, $\tau$. This is
trivially true if $i_k=j_k+1$ for some $k\leq r$ or $p_d=d+1$ for
some $d\leq r$, so we may assume that $i_l>j_l+1$ for all $1\leq
l\leq r$ and $p_l>l+1$ for all $l\leq r_0=\min\{r,s(\tau_r)\}$.

Suppose $\Supp{\sigma}\cap\Co_{r+1}=\emptyset$. If $\Supp{\tau}
\cap\Co_{r+1}=\emptyset$, too, then $\sigma_{r+1}$, $\tau_{r+1}$
satisfy~(\ref{formula:remained_case_ravno}). At the contrary, assume
$\Supp{\tau}\cap\Co_{r+1}=\{(p,r+1)\}\neq\emptyset$. We see that
$r_0=s(\tau_r)=d<r$ (if $r_0=r$, then
$(R_{\tau}^*)_{r+1,r}=r+1>r=(R_{\sigma}^*)_{r+1,r}$, so
$\sigma\ngtr^*\tau$, a contradiction). If
$\Supp{\tau}\cap\Ro_r=\emptyset$, then put
\begin{equation*}
\tau_0=(p_1,1)\ldots(p_d,d)(p,r)(p_{d+2},q_{d+2})\ldots(p_t,q_t),
\end{equation*}
i.e.,
$\Supp{\tau_0}=(\Supp{\tau}\setminus\{(p,r+1)\})\cup\{(p,r)\}$. In
this case, $\tau_0\in S_n^2$ and $\tau<^*\tau_0\leq^*\sigma$. If
$\tau_0=\sigma$, then $\tau=\sigma_{(p,r)}^{\to}\in N^0(\sigma)$. At
the same time, if $\tau_0<^*\sigma$, then $\tau\notin \wt L^0(\sigma)$.

On the other hand, assume
$\Supp{\tau}\cap\Ro_r=\{(r,j)\}\neq\emptyset$. If $i_j>r$, then
put $\tau_0$ to be the involution such that
\begin{equation*}
\Supp{\tau_0}=(\Supp{\tau}\setminus\{(r,j),(p,r+1)\})\cup\{(r+1,j),(p,r)\}.
\end{equation*}
Evidently, $\tau<^*\tau_0\leq^*\sigma$. If $\tau_0=\sigma$, then
$\tau=a\sigma_{(r,j)}^{(p,r+1)}\in N^0(\sigma)$, and if
$\tau_0<^*\sigma$, then $\tau\notin \wt L^0(\sigma)$.

Next, assume $\Supp{\tau}\cap\Ro_r=\{(r,j)\}\neq\emptyset$, but
$i_j=r$ (in particular, $\Supp{\sigma}\cap\Co_r=\emptyset$). It
follows from (\ref{formula:remained_case_ravno}d) that
$\Supp{\tau}\cap\Co_l\cap\Ro_i=\emptyset$ for all $j<l<r$,
$r<i\leq n$. Let $q=[(r-j+1)/2]$ and
\begin{equation*}
\begin{split}
z=\max\{&s\mid r-q+1\leq s<r\text{ and either
}\Supp{\tau}\cap\Ro_s=\emptyset\\
&\text{or }\Supp{\tau}\cap\Ro_s=\{(s, h)\},\text{ where }i_h>s\}.
\end{split}
\end{equation*}
It's easy to see that $z$ exists. Indeed, if
$\Supp{\sigma}\cap\Ro_s=\Supp{\tau}\cap\Ro_s\neq\emptyset$ for all
$r-q+1\leq s<r$, then
$\Supp{\sigma}\cap\Ro_{r-q+1}=\Supp{\tau}\cap\Ro_{r-q+1}=\{(r-q+1,r-q)\}$,
which contradicts our assumption. Further,
$\Supp{\tau}\cap\Co_z=\emptyset$, because $z\geq r-q+1>q\geq d$.
If $\Supp{\tau}\cap\Ro_z=\emptyset$ (resp.
$\Supp{\tau}\cap\Ro_z=\{(z,h)\}$), then put $\tau_0$ to be the
involution such that
\begin{equation*}
\begin{split}
&\Supp{\tau_0}=(\Supp{\tau}\setminus\{(p,r+1)\})\cup\{(p,z)\} \text{
(resp.}\\
&\Supp{\tau_0}=(\Supp{\tau}\setminus\{(z,h),(p,r+1)\})\cup\{(r+1,h),(p,z)\}).
\end{split}
\end{equation*}
Clearly, $\tau<^*\tau_0\leq^*\sigma$. If $\tau_0<^*\sigma$, then
$\tau\notin \wt L^0(\sigma)$. By the choice of $z$,
$\Supp{\sigma}\cap\Co_l=\emptyset$ for all $z<l\leq r+1$, so if
$\tau_0=\sigma$, then $\tau=\sigma_{(p,z)}^{\to}$ (resp.
$\tau=a\sigma_{(z,h)}^{(p,r+1)}$). In both cases, $\tau\in
N^0(\sigma)$.

\medskip iv) Now, let us consider the case when $\sigma_r$, $\tau_r$
satisfy (\ref{formula:remained_case_ravno}), $i_l>j_l+1$ for all
$1\leq l\leq r$, $p_l>l+1$ for all $l\leq r_0=\min\{r,s(\tau_r)\}$,
but $\Supp{\sigma}\cap\Co_{r+1}=\{(i,r+1)\}\neq\emptyset$. If
$\Supp{\tau}\cap\Co_{r+1}=\emptyset$, then $\sigma_{r+1}$,
$\tau_{r+1}$ satisfy (\ref{formula:remained_case_ravno}), so assume
$\Supp{\tau}\cap\Co_{r+1}=\{(p,r+1)\}\neq\emptyset$. Furthermore,
assume $r_0=s(\tau_r)=r$ and $p=p_{r+1}>p_r$ (and so $p<i$). Put
\begin{equation*}
\begin{split}
&\tau_0=(p_1,1)\ldots(p_{r-1},r-1)(p,r)(p_r,r+1)(p_{r+2},q_{r+2})\ldots(p_t,q_t),\text{
i.e.,}\\
&\Supp{\tau_0}=(\Supp{\tau}\setminus\{(p_r,r),(p,r+1)\})\cup\{(p,r),(p_r,r+1)\}.
\end{split}
\end{equation*}

Obviously, $\tau<^*\tau_0$. If $\tau_0<^*\sigma$, then $\tau\notin\wt L^0(\sigma)$, because $s(\tau_0)=s(\tau)=t=s(\sigma)$. Since
$p_r<p<i$, $\tau_0\neq\sigma$. If $\tau_0\nless^*\sigma$ (i.e.,
$i_r<p$), then put also
\begin{equation*}
\begin{split}
&\sigma_0=(i_1,j_1)\ldots(i_{r-1},r-1)(i,r)(i_r,r+1)(i_{r+2},j_{r+2})\ldots(i_t,j_t),\text{
i.e.,}\\
&\Supp{\sigma_0}=(\Supp{\sigma}\setminus\{(i_r,r),(i,r+1)\})\cup\{(i,r),(i_r,r+1)\}.
\end{split}
\end{equation*}
In this case, $\sigma_0\in S_n^2$, $\sigma_0>^*\sigma$,
$\sigma_0>^*\tau_0$ and $s(\sigma_0)=t=s(\tau_0)$, so, by the second
induction hypothesis, there exists $w_1\in N^0(\sigma_0)\cup
N^-(\sigma_0)$ such that $w_1\geq^*\tau_0$. One can check that this
implies an existence of $w\in N^0(\sigma)\cup N^-(\sigma)$ such that
$w\geq^*\tau$. (In fact, $w$ is obtained from $\sigma$ by
the ``same'' operation as $w_1$ from $\sigma_0$.)

\medskip v) Next, suppose $\sigma_r$, $\tau_r$
satisfy (\ref{formula:remained_case_ravno}), $i_l>j_l+1$ for all
$1\leq l\leq r$, $p_l>l+1$ for all $l\leq r_0=\min\{r,s(\tau_r)\}$,
$\Supp{\sigma}\cap\Co_{r+1}=\{(i,r+1)\}\neq\emptyset$,
$\Supp{\tau}\cap\Co_{r+1}=\{(p,r+1)\}\neq\emptyset$, but
$r_0=s(\tau_r)=d<r$. If $\Supp{\tau}\cap\Ro_r=\emptyset$, then put $\tau_0=(p_1,1)\ldots(p_d,d)(p,r)(p_{d+2},q_{d+2})\ldots(p_t,q_t)$, i.e.,
\begin{equation*}
\Supp{\tau_0}=(\Supp{\tau}\setminus\{(p,r+1)\})\cup\{(p,r)\}.
\end{equation*}
At the same time, if
$\Supp{\tau}\cap\Ro_r=\{(r,q)\}\neq\emptyset$, then define
$\tau_0$ to be the involution such that
\begin{equation*}
\Supp{\tau_0}=(\Supp{\tau}\setminus\{(r,q),(p,r+1)\})\cup\{(r+1,q),(p,r)\}.
\end{equation*}
In both cases, $\tau<^*\tau_0$ and $\tau_0\neq\sigma$. If
$\tau_0<^*\sigma$, then $\tau\notin\wt L^0(\sigma)$, because
$s(\tau_0)=s(\tau)=t=s(\sigma)$.

At the contrary, suppose $\tau_0\nless^*\sigma$. If
$\Supp{\sigma}\cap(\Co_r\cup\Ro_r)=\emptyset$, then define
$\sigma_0$ to be the involution such that
$\Supp{\sigma_0}=(\Supp{\sigma}\setminus\{(i,r+1)\})\cup\{(i,r)\}$.
If $\Supp{\sigma}\cap\Co_r=\{(x,r)\}\neq\emptyset$, then $x<r$, so
define $\sigma_0$ by putting
\begin{equation*}
\Supp{\sigma_0}=(\Supp{\sigma}\setminus\{(x,r),(i,r+1)\})\cup\{(i,r),(x,r+1)\}.
\end{equation*}
Finally, if $\Supp{\sigma}\cap\Ro_r=\{(r,y)\}\neq\emptyset$, then
define $\sigma_0$ by putting
\begin{equation*}
\Supp{\sigma_0}=(\Supp{\sigma}\setminus\{(r,y),(i,r+1)\})\cup\{(r+1,y),(i,r))\}.
\end{equation*}
In all cases, $\sigma_0\in S_n^2$, $\sigma_0>^*\sigma$,
$\sigma_0>^*\tau_0$ and $s(\sigma_0)=s(\tau_0)=t$. Hence, by the
second inductive assumption, there exists $w_1\in N^0(\sigma_0)\cup
N^-(\sigma_0)$ such that $w_1\geq^*\tau_0$. One can check that this
implies an existence of $w\in N^0(\sigma)\cup N^-(\sigma)$ such that
$w\geq^*\tau$. (In fact, $w$ is obtained from $\sigma$ by
the ``same'' operation as $w_1$ from $\sigma_0$.) The proof of
(\ref{formula:remained_ravno_r}) is complete.

\medskip vi) Now, we may assume without loss of generality that
$\sigma$, $\tau$ satisfy (\ref{formula:remained_case_ravno}). First,
suppose $i_k=j_k+1$ for some $k\leq t$ and $q_l>l+1$ for all $1\leq
l\leq d<k$. Then $(i_k,j_k)\in M(\sigma)$ and
$\tau<^*w=\sigma_{(i_k,j_k)}^-\in N^-(\sigma)$. Second, assume
$i_l>j_l+1$ for all $1\leq l\leq t$ and $p_l>q_l+1=l+1$ for all $l$ such that $q_l=l\leq j_t$. If
$q_t=j_t=t$, then $i_l\geq
p_l\geq p_t>t+1$ for all $1\leq l\leq t$, so
$\Supp{\sigma}\cap\Ro_{t+1}=\Supp{\tau}\cap\Ro_{t+1}=\emptyset$.
Thus, $\wt\tau=P_t(\tau)$ and $\wt\sigma=P_t(\sigma)$ belong to
$S_n^2\cap\wt S_{n-1}$, so, by the first induction hypothesis, there
exists $\wt w\in N^0(\wt\sigma)\cup N^-(\wt\sigma)$ such that
$\wt\tau\leq^*\wt w$. Let $w$ be the unique involution such that
$\wt w=P_t(w)$, then $w\in N^0(\sigma)\cup N^-(\sigma)$ and
$\tau\leq^*w$. On the other hand, if $q_t>j_t$, then, arguing as on step iii), one can show that either $\tau\in N^0(\sigma)$ or $\tau\notin\wt L^0(\sigma)$.

\medskip vii) Finally, let us consider the most interesting case when $\sigma$, $\tau$ satisfy (\ref{formula:remained_case_ravno}), but $p_d=d+1$
for some $d\leq t$ such that $q_d=d\leq j_t$. Put $\tau_0=\tau_{(d+1,d)}^-$, then
$\tau_0<^*\tau<^*\sigma$ and $s(\tau_0)=t-1<t=s(\sigma)$. By Lemma
\ref{lemm:hard_minus}, there exists $w_1\in N^-(\sigma)$ such that
$w_1\geq^*\tau_0$ (in fact, $w_1>^*\tau_0$). If $w_1>^*\tau$, then
the result follows. Thus, it remains to consider the case
$w_1\ngtr^*\tau$. This means that $(i_d,d)\in M(\sigma)$ and
$w_1=\sigma_{(i_d,d)}^-$. Denote $\sigma_0=w_1$ and
$\wt\sigma=P_{d-1}(\sigma_0)$, $\wt\tau=P_{d-1}(\tau_0)$. We see
that $\wt\sigma$, $\wt\tau\in\wt S_{n-1}\cap S_n^2$,
$\wt\sigma\geq^*\wt\tau$ and $s(\wt\sigma)=t-1=s(\wt\tau)$. If
$\wt\tau=\wt\sigma$, then $\Supp{\sigma}\cap\Co_{d+1}=\emptyset$
(if
$\Supp{\sigma}\cap\Co_{d+1}=\Supp{\tau}\cap\Co_{d+1}\neq\emptyset$,
then (\ref{formula:ortog_supp}) doesn't hold for $\tau$, because
$(d+1,d)\in\Supp{\tau}$). Moreover, by
(\ref{formula:remained_case_ravno}), $i_l\geq p_l>p_d=d+1$ for all
$1\leq l\leq d-1$, so $\Supp{\sigma}\cap\Ro_{d+1}=\emptyset$. In
other words, $\sigma(d+1)=d+1$; in particular,
\begin{equation*}
m=\max\{i\mid
i_d<i\leq d+1\text{ and }\sigma(i)=i\}
\end{equation*}
exists. Thus,
$w=\sigma_{(i_d,d)}^{\uparrow}\in N^0(\sigma)$ is well-defined and
$\tau\leq^*w$.

From now on, assume $\wt\tau<^*\wt\sigma$. Then the first
inductive assumption shows that there exists\linebreak $\wt w_2\in
N^0(\wt\sigma)\cup N^-(\wt\sigma)$ such that $\wt w_2\geq^*\wt\tau$.
Denote $i=i_d$. If $\Supp{\wt w_2}\cap(\Co_i\cup\Ro_i)=\emptyset$,
then denote by $w_2$ the unique involution such that $\wt
w_2=P_{d-1}(w_2)$ and define $w$ by putting
$\Supp{w}=\Supp{w_2}\cup\{(i,d)\}$. It's easy to see that $w\in
N^0(\sigma)\cup N^-(\sigma)$ and $w\geq^*\tau$. On the other hand,
suppose $(i,j)\in\Supp{\wt w_2}$ for some $j$. If $j<d$, then $\wt
w_2=\wt\sigma_{(x,j)}^{\uparrow}$ for some $x$. (Indeed, if $\wt
w_2=a\wt\sigma_{(x,j)}^{(\alpha,\beta)}$ for some $(\alpha,\beta)\in
A_{x,j}(\wt\sigma)$, then $\alpha=i$, which contradicts
(\ref{formula:ortog_supp}).) In this case, put
$w=b\sigma_{(i,d)}^{(\alpha,\beta)}$ for some $(\alpha,\beta)\in
B_{i,d}(\sigma)$ (since $x>i$, $(\alpha,\beta)$ exists). It follows
from $\wt w_2\geq^*\wt\tau$ that $w\geq^*\tau$.

Next, assume $(i,j)\in\Supp{\wt w_2}$ for some $j>d$. Denote
$(x,j)=\Supp{\sigma}\cap\Co_j$ (i.e., $\Supp{\wt
w_2}=(\Supp{\wt\sigma}\setminus\{(x,j)\})\cup\{i,j\}$). Since $\wt
w_2$ is well-defined, there are no $(\alpha,\beta)\in\Supp{\sigma}$
such that\linebreak $(\alpha,\beta)<(x,j)$ and $\alpha>i$. But
$(i,d)\in M(\sigma)$, so
$\Supp{\sigma}\cap\{(\alpha,\beta)\in\Phi\mid
(\alpha,\beta)<(i,d)\}=\emptyset$. Hence $(x,j)\in M(\sigma)$ and
$w=\sigma_{(x,j)}^-\geq^*\tau$ (in fact, $w>^*\tau$).

Similarly, assume $(r,i)\in\Supp{\wt w_2}$ for some $r>i$. Denote
$(r,s)=\Supp{\sigma}\cap\Ro_r$, i.e., $\Supp{\wt
w_2}=(\Supp{\wt\sigma}\setminus\{(r,s)\})\cup\{r,i\}$ (it follows
from (\ref{formula:remained_case_ravno}) that $s>d$). If
$\sigma(a)=a$ for some $i<a<d$, then $w=\sigma_{(i,d)}^{\uparrow}$
is well-defined and $w\geq^*\tau$, so suppose $\sigma(a)\neq a$ for
all $i<a<d$. Since $\wt w_2$ is well-defined and $(i,d)\in
M(\sigma)$, there are no $(\alpha,\beta)\in\Supp{\sigma}$ such that
$(\alpha,\beta)<(i,d)$ or $(\alpha,\beta)<(r,s)$,
$(\alpha,\beta)\nless(r,i)$. Thus, $w=a\sigma_{(i,d)}^{(r,s)}$ is
well-defined and $w\geq^*\tau$. The proof is complete.}

\medskip\lemmp{Let $\sigma\in S_n^2$. Then $L^+(\sigma)=N^+(\sigma)$.
\label{lemm:hard_plus}}{By
Lemma \ref{lemm:N_plus_in_L_plus}, it's enough to check that
$L^+(\sigma)\subseteq N^+(\sigma)$. By definition,
\begin{equation*}
L^+(\sigma)=\{\sigma'\in S_n^2\mid \sigma'\leq^*\sigma,\text{
}s(\sigma')>s(\sigma),\text{ and if }\sigma'\leq^* w<^*\sigma,\text{
then }w=\sigma'\},\\
\end{equation*}
so it suffice to show that
\begin{equation}
\begin{split}
&\text{if }
\tau\leq^*\sigma\text{ and }s(\tau)>s(\sigma),\\
&\text{then there exists }\sigma'\in N^+(\sigma)\cup N^0(\sigma)\cup N^-(\sigma)
\text{ such that }\tau\leq^*\sigma'<^*\sigma.
\end{split}\label{formula:proof_plus}
\end{equation}
We will proceed by induction on $n$
(for $n=1$, there is nothing to prove). The proof is rather long, so we split it into
four steps.

\medskip i) Note that if $\sigma=w_0$, the maximal element of $S_n^2$ with respect to $\leq^*$, then there is nothing to prove, because $L^+(w_0)=\emptyset$. Hence we may use the second (downward) induction on $\leq^*$. Let $\sigma=(i_1,j_1)\ldots(i_s,j_s)\in S_n^2$, $\sigma<^*w_0$,
$\tau=(p_1,q_1)\ldots(p_t,q_t)<^*\sigma$ and $s<t$.
Consider the following conditions (cf. (\ref{formula:remained_case_minus}) and (\ref{formula:remained_case_ravno})):
\begin{equation}
\begin{split}
&\text{a) There exists }k\leq s\text{ such that }
i_l>l+1\text{ for all }1\leq l\leq k-1.\\
&\text{b) There exists }d\leq k\text{ such that }q_l=l\text{ for all }
1\leq l\leq d\\
&\hphantom{\text{b) }}\text{ and either }d<k,\text{ }q_{d+1}>k,
\text{ or }d=k,\text{ }p_d=d+1.\\
&\text{c) }j_l=l\text{ for all }1\leq l\leq d,\text{ and if }
d<k,\text{ then either }i_k=j_k+1\text{ or }k=s.\\
&\text{d) }(p_l,q_l)=(p_l,l)>(p_{l+1},q_{l+1})=
(p_{l+1},l+1),\text{ i.e., }p_l>p_{l+1},\text{ for all }1\leq l\leq d-1.\\
&\text{e) }(i_l,j_l)=(i_l,l)\geq(p_l,q_l)=
(p_l,l),\text{ i.e., }i_l\geq p_l,\text{ for all }1\leq l\leq d.
\end{split}\label{formula:remained_case_plus}
\end{equation}
Pick $r\leq s$. Define $\sigma_r$ and $\tau_r$ as in the previous Lemmas. We claim that
\begin{equation}
\begin{split}
&\text{if (\ref{formula:proof_plus}) holds for all }
\sigma,\tau\text{ satisfying (\ref{formula:remained_case_plus})},\\
&\text{then (\ref{formula:proof_plus})
holds for all }\sigma,\tau\in S_n^2.
\end{split}\label{formula:remained_plus_r}
\end{equation}
Clearly,
it's enough to prove that if (\ref{formula:proof_plus}) holds for all
$\sigma $, $\tau$ satisfying (\ref{formula:remained_case_plus}), and $\sigma_r$, $\tau_r$
don't satisfy~(\ref{formula:remained_case_plus}) for some $1\leq r\leq s$,
then (\ref{formula:proof_plus}) holds for $\sigma$,
$\tau$.

We will proceed by induction on $r$. Evidently,
there exist $w_1=\sigma$, $w_2$, $\ldots$, $w_z=\tau\in S_n^2$
such that $w_1>^*w_2>^*\ldots>^*w_z$ and $w_{i+1}\in L^+(w_i)$
for all $1\leq i<z$, so we may assume $\tau\in L^+(\sigma)$. The base $r=1$
is evident (see step ii) of the proof of Lemma \ref{lemm:hard_minus}). The induction step can be performed as on steps iii)--v) of the proof of Lemma \ref{lemm:hard_ravno}, so assume without loss of generality that
$\sigma$, $\tau$ satisfy (\ref{formula:remained_case_plus}).

\medskip ii) If either $i_k=j_k+1$ for some $k\leq t$ and $q_l>l+1$ for all $1\leq
l\leq d<k$, or
$i_l>j_l+1$ for all $1\leq l\leq s$ and $p_l>q_l+1=l+1$ for all $l$ such that $q_l=l\leq j_s$, then, arguing as on step vi) of the proof of Lemma \ref{lemm:hard_ravno}, we obtain the result. Let us consider the most interesting case when $\sigma$, $\tau$ satisfy~(\ref{formula:remained_case_plus}), but $p_d=d+1$
for some $d\leq t$ such that $q_d=d\leq j_s$ (and so $j_d=d$, too).

Put $\tau_0=\tau_{(d+1,d)}^-$, then
$\tau_0<^*\tau<^*\sigma$ and $s(\tau_0)=t-1$. First, suppose $s(\sigma)=s=t-1=s(\tau_0)$. By Lemma
\ref{lemm:hard_ravno}, there exists $w_1\in N^0(\sigma)\cup N^-(\sigma)$ such that
$w_1\geq^*\tau_0$. If $w_1>^*\tau$, then
the result follows. Thus, it remains to consider the case
$w_1\ngtr^*\tau$ (clearly, $w_1\neq\tau$). Assume $w_1\in N^0(\sigma)$, then either $w_1=\sigma_{(i_d,d)}^{\to}$ or $w_1=a\sigma_{(x,y)}^{(i_d,d)}$ for some $(x,y)\in\Supp{\sigma}$, $y<d$. If $w_1=\sigma_{(i_d,d)}^{\to}$, then there exists $m>d$ such that $\sigma(m)=m$ and $\sigma(l)\neq l$ for all $d<l<m$. Suppose $l>d+1$, then there exists $j+1\leq p<i_d$ such that $(p, m)\in C_{i,j}(\sigma)$, and so $w=c\sigma_{(i,j)}^{p,m}>^*\tau$.

At the contrary, suppose $m=d+1$ and $\sigma(l)\neq l$ for all $d+1<l<i$. Suppose
\begin{equation*}
r=\min\{l\mid d+1<l<i\text{ and }(z,l)\in\Supp{\sigma}\text{ for some }z<i\}
\end{equation*}
exists (so $(z,r)\in B_{i,j}(\sigma)$). We claim that $w=b\sigma_{(z,l)}^{(i,j)}>^*\tau$. Indeed, it suffice to show that\linebreak$(R_w^*)_{i,j}\geq(R_{\tau}^*)_{i,j}$ for all $z<i<i_d$, $d<j<r$. Denote by $\gamma$ (resp. $\gamma'$) the number of $(a,b)\in\Supp{\sigma}$ (resp. $(a,b)\in\Supp{\tau}$) such that $b<d$ and $a<i$. Since $\sigma(l)\neq l$ for all $j+1<l\leq j$, there are $((j-d-1)-\gamma)$ (resp. at most $((j-d-1)-\gamma')$) elements $(a,b)$ in $\Supp{\sigma}$ (resp. in $\Supp{\tau}$) such that $d+1<b\leq j$ and $a>i_d$ (resp. $d+1<b\leq j$ and $a\geq i$). By (\ref{formula:remained_case_plus}e), $\gamma\leq\gamma'$. Thus,
\begin{equation*}
\begin{split}
(R_w^*)_{i,j}-(R_{\tau}^*)_{i,j}&=\#\{(a,b)\in\Supp{\sigma}\mid b<d\text{ and }a>i\}\\
&-\#\{(a,b)\in\Supp{\tau}\mid b>d+1\text{ and }a\geq i\}\\
&\geq((j-d-1)-\gamma)-((j-d-1)-\gamma')=\gamma'-\gamma\geq0,
\end{split}
\end{equation*}
as required. If $r$ doesn't exist, then, arguing as above, one can check that $w=\sigma_{(i_d,d)}^{\uparrow}>^*\tau$. (In fact,
$\Supp{w}=(\Supp{\sigma}\setminus\{(i_d,d)\})\cup\{(d+1,d)\}$.) The case $w_1=a\sigma_{(x,y)}^{(i_d,d)}$ is similar.

\medskip iii) Next, suppose $w_1\in N^-(\sigma)$. This means that $(i_d,d)\in M(\sigma)$ and $\Supp{w_1}=\Supp{\sigma}\setminus\{(i_d,d)\}$. Denote $\sigma_0=w_1$, $\tau_0=\tau_{(d+1,d)}^-$ and
$\wt\sigma=P_{d-1}(\sigma_0)$, $\wt\tau=P_{d-1}(\tau_0)$. We see
that $\wt\sigma$, $\wt\tau\in\wt S_{n-1}\cap S_n^2$,
$\wt\sigma>^*\wt\tau$ and $s(\wt\sigma)=s-1<t-1=s(\wt\tau)$. The first
inductive assumption shows that there exists $\wt w_2\in
N^+(\wt\sigma)\cup N^0(\wt\sigma)\cup N^-(\wt\sigma)$ such that $\wt w_2\geq^*\wt\tau$. If $\wt w_2\in N^-(\wt\sigma)\cup N^0(\wt\sigma)$, then there exists $w\in N^-(\sigma)\cup N^0(\sigma)$ such that $w\geq^*\tau$, see step vi) of the proof of Lemma \ref{lemm:hard_minus} and step vii) of the proof of Lemma \ref{lemm:hard_ravno}, so it remains to consider the case $\wt w_2\in N^+(\wt\sigma)$.

If $\wt w_2=c\wt\sigma_{(i,j+1)}^{\alpha,\beta}$ for some $(i,j)\in\Supp{\sigma}$, $(i,j)>(i_d,d)$, but $(\alpha,j)\ngtr(i_d,d)$, $(i,\beta)\ngtr(i_d,d)$, then $d+1<\alpha<i_d$ and $\beta\leq d<i_d$, so $\sigma(l)\neq l$ for all $\alpha<l<i_d$. Hence $w=\sigma_{(i_d,d)}^{\uparrow}\geq^*\tau$ (in fact, $\Supp{w}=(\Supp{\sigma}\setminus\{(i_d,d)\})\cup\{(\alpha,d)\}$). On the other hand, if $\wt w_2=c\wt\sigma_{(i,j)}^{i_d,\beta}$ for some $(i,j)\in\Supp{\sigma}$, then $i>i_d$ and $j>d$, so $w=\sigma_{(i,j)}^{\to}\geq^*\tau$ (in fact, $\Supp{w}=(\Supp{\sigma}\setminus\{(i,j)\})\cup\{(i,\beta)\}$). If $\wt w_2=c\wt\sigma_{(i,j)}^{\alpha,i_d}$ for some $(i,j)\in\Supp{\sigma}$, then $i>i_d$, $j>d$ and one can easily check that $(i,j)\in A_{(i_d,d)}(\sigma)$, so $w=a\sigma_{(i_d,d)}^{(i,j)}\geq^*\tau$. In all other cases, $w\geq^*\tau$, where $\Supp{w}=\Supp{w_2}\cup\{(i_d,d)\}$ and $w_2$ is the unique involution such that $\wt w_2=P_{d-1}(w_2)$.

\medskip iv) Finally, assume $s(\sigma)=s<t-1=s(\tau_0)$ (see the beginning of step ii)). Suppose
\begin{equation*}
r=\min\{l\mid(z,l)\in\Supp{\tau},\text{ }l>d+1\text{ and }z\leq i_d\}
\end{equation*}
exists. Then $\sigma\geq^*\tau_1>^*\tau$, where
\begin{equation*}
\Supp{\tau_1}=(\Supp{\tau}\setminus\{(d,d+1),(z,r)\})\cup\{(z,d)\},
\end{equation*}
so $\tau\notin L^+(\sigma)$. At the same time, if $r$ doesn't exists, then define $\sigma_0$ to be the involution such that $\Supp{\sigma_0}=\Supp{\sigma}\setminus\{(i_d,d)\}$. In this case, $\sigma_0\geq^*\tau_0$ and $s(\sigma_0)=s-1<t-1=s(\tau_0)$, so we may proceed as on the previous step. The proof is complete.}

\medskip Finally, we will prove the fact used in the proof of Theorem \ref{mcoro}. \lemmp{Let $\sigma\in S_n^2$. Then $L'(\sigma)=N'(\sigma)$.\label{lemm:hard_bis}}{Pick an involution $\tau=\sigma{(i,j)}^-\in N^-(\sigma)=L^-(\sigma)$. If $\tau\notin N'(\sigma)$, then there exists $m$ such that $j\leq m\leq i$ and $\sigma(m)=m$. Assume, for example, that $i\neq m$, then $\tau<^*w<^*\sigma$, where $$\Supp{w}=\Supp{\tau}\cup\{(i,m)\},$$ so $\tau\notin L'(\sigma)$. The case $j\neq m$ is similar.

On the other hand, suppose $\tau\notin L'(\sigma)$. Then there exists $w\in L^+(\sigma)\cup L^0(\sigma)$ such that $\tau<^* w<^*\sigma$. By Lemmas \ref{lemm:hard_plus} and \ref{lemm:hard_ravno}, $L^+(\sigma)=N^+(\sigma)$ and $L^0(\sigma)=N^0(\sigma)$ respectively. By (\ref{formula:ravno_mezhdu}),
\begin{equation*}
\Supp{\tau}\cap Y=\Supp{w}\cap Y=\Supp{\sigma}\cap Y,
\end{equation*}
where $Y=\{(p,q)\in\Phi\mid(p,q)\nleq(i,j)\}$, so $\Supp{w}\cap\wt Y\neq\emptyset$, where $\wt Y=\{(p,q)\in\Phi\mid(p,q)\leq(i,j)\}$. It follows from this fact and the definitions of $N^+(\sigma)$ and $N^0(\sigma)$ (see Subsection \ref{sst:Near}) that there exists $m$ such that $j\leq m\leq i$ and $\sigma(m)=m$, so $\tau\notin N'(\sigma)$. (For example, if $w=\sigma_{(\alpha,\beta)}^{\uparrow}$ for some $(\alpha,\beta)\in\Supp{\sigma}$, then $\Supp{w}\setminus\Supp{\sigma}=\{(m,\beta)\}$.) This concludes the proof.}

\sect{Concluding remarks}\label{sect:remarks}
\sst\label{sst:dim_orbit} Let $\sigma\in S_n^2$. Using results of \cite{Panov}, one can easily obtain a formula for the dimension of the orbit $\Omega_{\sigma}$. Let $\xi\colon\Supp{\sigma}\to\Cp^{\times}$ be a map. Recall the notation from Subsection~\ref{sst:repr_theo_approach}. As above, let $l(\sigma)$ be the length of a reduced expression of $\sigma$ as a product of simple reflections, and $s(\sigma)=|\Supp{\sigma}|$ (obviously, if $\Supp{\sigma}=\{(i_1,j_1),\ldots,(i_t,j_t)\}$, then $s(\sigma)=t$). By \cite[Theorem 1.2]{Panov}, $\Theta_{\sigma,\xi}$ is an irreducible affine variety of dimension $\dim\Theta_{\sigma,\xi}=l(\sigma)-s(\sigma)$. By Lemma \ref{lemm:unip_orbits}, $\Omega_{\sigma}=\bigcup_{\xi}\Theta_{\sigma,\xi}$. Denote $\Theta_0=\Theta_{\sigma,\xi_0}$, where $\xi_0(\alpha)=1$ for all $\alpha\in\Supp{\sigma}$ (in other words, $\Theta_0$ is the $U$-orbit of $X_{\sigma}^t$).

\medskip\propp{Let $\sigma\in S_n$ be an involution. Then $\dim\Omega_{\sigma}=l(\sigma)$.\label{prop:dim_Omega}}{Let $Z=\mathrm{Stab}_BX_{\sigma}^t$ be the stabilizer of $X_{\sigma}^t$ in $B$. One has
\begin{equation*}
\dim\Omega_{\sigma}=\dim B-\dim Z.
\end{equation*}
Recall that $B=U\rtimes D$. Suppose $g=ud\in Z$, where $u\in U$, $d\in D$, then $g.X_{\sigma}^t=u.(d.X_{\sigma}^t)=X_{\sigma}^t$. But $d.X_{\sigma}^t=f_{\sigma,\xi}$, where $\xi_l=d_{i_l,i_l}/d_{j_l,j_l}$ (cf. the proof of Lemma \ref{lemm:unip_orbits}), so $g.X_{\sigma}^t=u.f_{\sigma,\xi}\in\Theta_{\sigma,\xi}$. Since $\Theta_{\sigma,\xi}\neq\Theta_0$ for $\xi\neq\xi_0$, we conclude that $\xi=\xi_0$, so $d.X_{\sigma}^t=X_{\sigma}^t$. Hence $d\in Z_D$ and $u\in Z_U$, where $Z_D=\mathrm{Stab}_DX_{\sigma}^t$ (resp. $Z_U=\mathrm{Stab}_UX_{\sigma}^t$) is the stabilizer of $X_{\sigma}^t$ in $D$ (resp. in $U$).

Since $B=U\rtimes D$ as algebraic groups, the maps
\begin{equation*}
\begin{split}
&\phi\colon B\to U\times D\colon g=ud\mapsto(u,d)\text{ and}\\
&\psi\colon U\times D\to B\colon(u,d)\mapsto ud
\end{split}
\end{equation*}
are inverse isomorphisms of algebraic varieties. We checked that $\phi(Z)\subseteq Z_U\times Z_D$. The inclusion $\psi(Z_U\times Z_D)\subseteq Z$ is evident, so $\phi(Z)=Z_U\times Z_D$. Thus, $\dim Z=\dim Z_U+\dim Z_D$.

But $d\in D$ belongs to $Z_D$ if and only if $\xi=\xi_0$, i.e., $d_{i_l,i_l}=d_{j_l,j_l}$ for all $1\leq l\leq t$. Hence $$\dim Z_D=\dim D-|\Supp{\sigma}|=n-s(\sigma).$$
On the other hand, since $\dim\Theta_0=l(\sigma)-s(\sigma)$, we obtain
$$\dim Z_U=\dim U-\dim\Theta_0=\dim B-n-l(\sigma)+s(\sigma).$$ We conclude that $\dim Z=\dim B-l(\sigma)$, and so $\dim\Omega_{\sigma}=l(\sigma)$, as required.}

\sst\label{sst:closure_conj} Here we present a conjectural description of the closure of a given $B$-orbit $\Omega_{\sigma}$, $\sigma\in S_n^2$. Namely, we describe a subvariety $Z_{\sigma}\subseteq\nt^*$, prove that $\overline{\Omega}_{\sigma}\subseteq Z_{\sigma}$ and prove the equality $\overline{\Omega}_{\sigma}=Z_{\sigma}$ in some particular cases. Define $S_{\sigma}$ to be the set of maximal elements from $\Supp{\sigma}$ with respect to the order $\leq$ on $\Phi$. Let $$\Mo_{\sigma}=\{\alpha\in\Phi\mid\alpha>\beta\text{ for some }\beta\in S_{\sigma}\}.$$

\exam{i) Let $\sigma=w_0=(n,1)(n-1,2)\ldots(n-n_0+1,n_0)$, $n_0=[n/2]$, be the maximal element of $S_n^2$ with respect to $\leq_B$. Then $S_{\sigma}=\{(n,1)\}$ and $\Mo_{\sigma}=\emptyset$.

ii) Let $n=8$, $\sigma=(5,1)(7,3)(6,4)$. Then $S_{\sigma}=\{(5,1), (7,3)\}$ and $\Mo_{\sigma}=\{(6,1)$, $(7,1)$, $(8,1)$, $(7,2)$, $(8,2)$, $(8,3)\}$ (these elements are grey on the picture below).
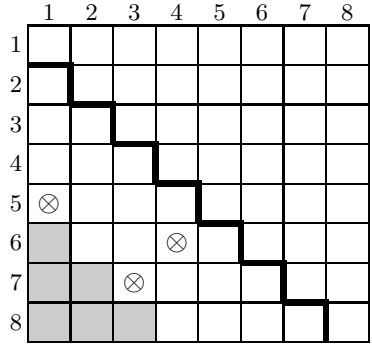
\begin{equation*}
\mymatrix{\pho& \pho& \pho& \pho& \pho& \pho& \pho& \pho\\
\Top{2pt}\Rt{2pt}\pho& \pho& \pho& \pho& \pho& \pho& \pho& \pho\\
\pho& \Top{2pt}\Rt{2pt}\pho& \pho& \pho& \pho& \pho& \pho& \pho\\
\pho& \pho& \Top{2pt}\Rt{2pt}\pho& \pho& \pho& \pho& \pho& \pho\\
\otimes& \pho& \pho& \Top{2pt}\Rt{2pt}\pho& \pho& \pho& \pho& \pho\\
\gray\pho& \pho& \pho& \otimes& \Top{2pt}\Rt{2pt}\pho& \pho& \pho& \pho\\
\gray\pho& \gray\pho& \otimes& \pho& \pho& \Top{2pt}\Rt{2pt}\pho& \pho& \pho\\
\gray\pho& \gray\pho& \gray\pho& \pho& \pho& \pho& \Top{2pt}\Rt{2pt}\pho& \pho\\
}\end{equation*}}

Let $A\in\nt^*$. To each $(r,s)\in\Phi$ one can assign polynomial $\gamma_{r,s}$ in $A_{i,j}$ of the form
\begin{equation*}
\gamma_{r,s}(A)=(A^2)_{r,s}=\sum_{k=1}^nA_{r,k}A_{k,s}=\sum_{k=s+1}^{r-1}A_{r,k}A_{k,s}.
\end{equation*}

Denote by $Z_{\sigma}$ the subvariety of $\nt^*$ defined by
\begin{equation}
\begin{split}
&\rk\pi_{i,j}(A)\leq(R_{\sigma}^*)_{i,j}\text{ for all }(i,j)\in\Phi,\\
&\gamma_{i,j}(A)=0\text{ for all }(i,j)\in\Mo_{\sigma}.
\end{split}\label{formula:conj_closure}
\end{equation}

\propp{Let $\sigma\in S_n$ be an involution. Then $\overline{\Omega}_{\sigma}\subseteq Z_{\sigma}$.\label{prop:Omega_in_Z}}{Let $A\in\overline{\Omega}_{\sigma}$. Lemma \ref{lemm:rank_through_orbit} guarantees that $A$ satisfies $\rk\pi_{i,j}(A)\leq(R_{\sigma}^*)_{i,j}$ for all $(i,j)\in\Phi$, so it remains to check that $\gamma_{i,j}(A)=0$ for all $(i,j)\in\Mo_{\sigma}$. Pick an element $(r,s)\in\Mo_{\sigma}$. Suppose $A\in\Omega_{\sigma}$. Recall that there exists $g\in B$ such that $A=(gX_{\sigma}^tg^{-1})_{\mathrm{low}}$, the strictly lower-triangular part of $y=gX_{\sigma}^tg^{-1}$ (see Subsection \ref{sst:our_order}). But $(X_{\sigma}^*)^2=0$, so $y^2=0$. In particular, $$(y^2)_{r,s}=\sum_{k=1}^ny_{r,k}y_{k,s}=\sum_{k=1}^sy_{r,k}y_{k,s}+\gamma_{r,s}(A)+\sum_{k=r}^ny_{r,k}y_{k,s}=0.$$

Clearly, $(R_{\sigma}^*)_{i,j}=0$ for all $(i,j)\in\Mo_{\sigma}$. Hence $y_{r,k}=A_{r,k}=0$ for all $1\leq k\leq s$, because if $1\leq k\leq s$, then $(r,k)>(r,s)\in\Mo_{\sigma}$, and so $(r,k)\in\Mo_{\sigma}$. This implies $\sum_{k=1}^sy_{r,k}y_{k,s}=0$. Similarly, $y_{k,s}=A_{k,s}=0$ for all $r\leq k\leq n$, because if $r\leq k\leq n$, then $(k,s)>(r,s)\in\Mo_{\sigma}$, and so $(k,s)\in\Mo_{\sigma}$. This implies $\sum_{k=r}^ny_{r,k}y_{k,s}=0$. Thus, $\gamma_{r,s}(A)=0$ for all $A\in\Omega_{\sigma}$, and so for all $A\in\overline{\Omega}_{\sigma}$.}

\medskip\hypo{Let $\sigma\in S_n$ be an involution. Then $\overline{\Omega}_{\sigma}=Z_{\sigma}$.\label{conj_closure}}

\medskip\nota{Suppose $\tau\leq^*\sigma$. Then $\Mo_{\sigma}\subseteq\Mo_{\tau}$, so that one can see immediately $\gamma_{r,s}(A)=0$ for $(r,s)\in\Mo_{\sigma}$ and $A\in\Omega_{\tau}$.}

\medskip Unfortunately, we can neither prove the irreducibility of $Z_{\sigma}$ nor compute its dimension. On the other hand, in some particular cases the proof of the equality $\overline{\Omega}_{\sigma}=Z_{\sigma}$ is more or less straightforward. Namely, assume $\Supp{\sigma}=\{(i_1,j_1),\ldots,(i_t,j_t)\}$ is a \emph{chain}, i.e., $(i_1,j_1)>\ldots>(i_t,j_t)$. (For instance, $\sigma=w_0$, or, more generally, $\sigma=(n,1)(n-1,2)\ldots(n-k+1,k)$ is maximal among all involutions with $k\leq n_0=[n/2]$ disjoint cycles.)
\propp{If $\Supp{\sigma}$ is a chain, then $\overline{\Omega}_{\sigma}=Z_{\sigma}$.\label{prop:conj_closure}}{In this case, $S_{\sigma}=\{(i_1,j_1)\}$, so $$\Mo_{\sigma}=\{(i,j)\in\Phi\mid i\geq i_1\text{ and }j\leq j_1\}\setminus\{(i_1,j_1)\}.$$ Suppose $A\in Z_{\sigma}$. Obviously, if $i>i_1$ or $j<j_1$, then $(R_{\sigma}^*)_{i,j}=0$, hence $\rk\pi_{i,j}(A)=0$ and so $A_{i,j}=0$. It follows that if $(r,s)\in\Mo_{\sigma}$, then $\gamma_{r,s}(A)=0$. Thus, $A\in\nt^*$ belongs to $Z_{\sigma}$ if and only if $\rk\pi_{i,j}(A)\leq(R_{\sigma}^*)_{i,j}$ for all $(i,j)\in\Phi$.

We need some more notation. For $1\leq i,j\leq n$, let $\wh\pi_{i,j}\colon\gt\to\gt$ be the map sending a matrix\linebreak $y\in\gt=\Mat_n(\Cp)$ to its upper-left $i\times j$ submatrix. Denote also by $P\colon\gt\to\gt$ the map defined by $P(y)_{i,j}=y_{n-j+1,i}$, $y\in\gt$, $1\leq i,j\leq n$. (Note that $\wh\pi_{i,j}=\pi_{i,j}\circ P$ for $(i,j)\in\Phi$.) Put $w=w_0\sigma$ and
\begin{equation*}
\begin{split}
&\Do(w)=\{(i,j)\mid w(i)>j\text{ and }w^{-1}(j)>i\},\\
&\Eo(w)=\{(i,j)\in\Do(w)\mid(i+1,j)\notin\Do(w)\text{ and }(i,j+1)\notin\Do(w)\},\\
&Z'=\{y\in\gt\mid\rk\wh\pi_{i,j}(y)\leq\rk\wh\pi_{i,j}(\dot w)\text{ for all }1\leq i,j\leq n\},\\
&Z=Z'\cap P(\nt^*)=\{y\in Z'\mid y_{i,j}=0\text{ for all }i\geq n-j+1\},\\
&Z''=\{y\in\gt\mid\rk\wh\pi_{i,j}(y)\leq\rk\wh\pi_{i,j}(\dot w)\text{ for all }(i,j)\in\Eo(w)\}.\\
\end{split}
\end{equation*}
(Clearly, $Z_{\sigma}=P(Z)$, because $\rk\wh\pi_{i,j}(\dot w)=\rk\pi_{i,j}(P(\dot w))=(R_{\sigma}^*)_{n-j+1,i}$.) For example, if $n=8$, $\sigma=(8,2)(6,3)$, then $w=(8,7,2,1)(5,4)$ (here we write $w$ as a product of disjoint cycles),
\begin{equation*}
\begin{split}
\Do(w)=\{&(1,2),(1,3),(1,4),(1,5),(1,6),\\
&(1,7),(3,2),(4,2),(4,4),(5,2),(6,2)\},\\
\Eo(w)=\{&(1,7),(4,4),(6,2)\}.
\end{split}
\end{equation*}
On the left (resp. right) picture below we draw $X_{\sigma}^*$ (resp. $\dot w$) as a rook placement. On the right picture boxes from $\Do(w)$ are grey and boxes from $\Eo(w)$ are marked by $\times$'s.
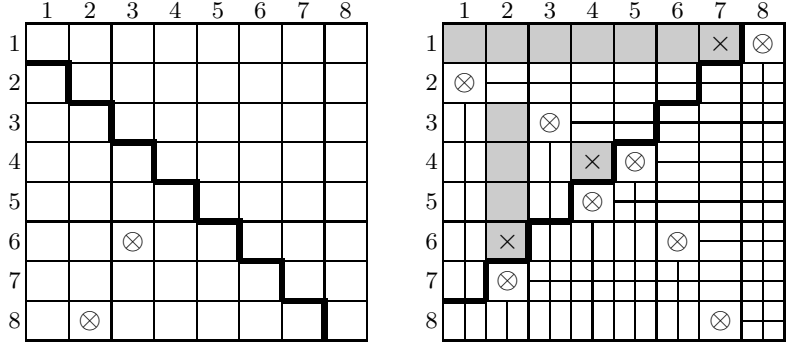
\begin{equation*}
\mymatrix{\pho& \pho& \pho& \pho& \pho& \pho& \pho& \pho\\
\Top{2pt}\Rt{2pt}\pho& \pho& \pho& \pho& \pho& \pho& \pho& \pho\\
\pho& \Top{2pt}\Rt{2pt}\pho& \pho& \pho& \pho& \pho& \pho& \pho\\
\pho& \pho& \Top{2pt}\Rt{2pt}\pho& \pho& \pho& \pho& \pho& \pho\\
\pho& \pho& \pho& \Top{2pt}\Rt{2pt}\pho& \pho& \pho& \pho& \pho\\
\pho& \pho& \otimes& \pho& \Top{2pt}\Rt{2pt}\pho& \pho& \pho& \pho\\
\pho& \pho& \pho& \pho& \pho& \Top{2pt}\Rt{2pt}\pho& \pho& \pho\\
\pho& \otimes& \pho& \pho& \pho& \pho& \Top{2pt}\Rt{2pt}\pho& \pho\\
}\qquad
\mymatrix{\gray\pho& \gray\pho& \gray\pho& \gray\pho& \gray\pho& \gray\pho& \Rt{2pt}\gray\times& \otimes\\
\otimes& \hpal\pho& \hpal\pho& \hpal\pho& \hpal\pho& \Rt{2pt}\hpal\pho& \Top{2pt}\hpal\pho& \hpal\vpal\pho\\
\vpal\pho& \gray\pho& \otimes& \hpal\pho& \Rt{2pt}\hpal\pho& \Top{2pt}\hpal\pho& \hpal\pho& \hpal\vpal\pho\\
\vpal\pho& \gray\pho& \vpal\pho& \Rt{2pt}\gray\times& \Top{2pt}\otimes& \hpal\pho& \hpal\pho& \hpal\vpal\pho\\
\vpal\pho& \gray\pho& \Rt{2pt}\vpal\pho& \Top{2pt}\otimes& \hpal\vpal\pho& \hpal\pho& \hpal\pho& \hpal\vpal\pho\\
\vpal\pho& \Rt{2pt}\gray\times& \Top{2pt}\vpal\pho& \vpal\pho& \vpal\pho& \otimes& \hpal\pho& \hpal\vpal\pho\\
\Rt{2pt}\vpal\pho& \Top{2pt}\otimes& \hpal\vpal\pho& \hpal\vpal\pho& \hpal\vpal\pho& \hpal\vpal\pho& \hpal\pho& \hpal\vpal\pho\\
\Top{2pt}\vpal\pho& \vpal\pho& \vpal\pho& \vpal\pho& \vpal\pho& \vpal\pho& \otimes& \hpal\vpal\pho\\
}\end{equation*}

Evidently, $Z'\subseteq Z''$. (In \cite{MillerSturmfels}, $Z'$ is called a \emph{determinantal matrix Schubert variety}.) But it follows from \cite[Theorem 15.15]{MillerSturmfels} that $Z'=Z''$. Furthermore, \cite[Theorem 15.31]{MillerSturmfels} claims that $Z'$ is a smooth irreducible affine subvariety of $\gt$ of dimension $\dim Z'=n^2-l(w)$. Denote $$V=\{y\in\gt\mid y_{i,j}=0\text{ for all }i<n-j+1\}.$$ Since $\Supp{\sigma}$ is a chain, $\Do(w)$ is contained in the set $\{(i,j)\mid i<n-j+1\}$. This means that $Z'\cong Z\times V$ as affine varieties, hence $Z$ is a smooth irreducible affine variety of dimension $$\dim Z=\dim Z'-\dim V=\dim\nt^*-l(w).$$ Thus, $Z_{\sigma}=P(Z)$ is an irreducible subvariety of $\nt^*$ of dimension $\dim Z=\dim\nt^*-l(w)$. But it is well-known that $l(w)=l(w_0\sigma)=|\Phi|-l(\sigma)=\dim\nt^*-l(\sigma)$, so $\dim Z_{\sigma}=l(\sigma)$. Finally, $\overline{\Omega}_{\sigma}\subseteq Z_{\sigma}$ (by Proposition \ref{prop:Omega_in_Z}) and $\dim\overline{\Omega}_{\sigma}=l(\sigma)$ (by Proposition \ref{prop:dim_Omega}), so $\overline{\Omega}_{\sigma}=Z_{\sigma}$.}

\medskip Note that Conjecture \ref{conj_closure} together with Conjecture \ref{conj_cones} imply an explicit description of the tangent cone $C_{\sigma}$ to the Schubert variety $X_{\sigma}$ for an involution $\sigma$.

\end{document}